\documentclass[11pt, oneside, letterpaper, fleqn, leqno]{article}

\usepackage[utf8]{inputenc}
\usepackage[T1]{fontenc}
\usepackage[top=1in, bottom=1in, left=1in, right=1in]{geometry}
\usepackage{fouriernc}
\usepackage{authblk}
\usepackage{mathtools}
\usepackage{multicol}
\usepackage{graphicx}
\usepackage{amssymb}
\usepackage{bbold}
\usepackage{appendix}

\newcommand{\ket}[1]{\rvert #1 \rangle}

\newcommand{\kket}[2]{\rvert #1 \rangle_{#2}}

\newcommand{\bbraket}[4]{{}_{#1} \langle #2 \rvert #3 \rangle_{#4}}
\newcommand{\BBraket}[5]{{}_{#1} \langle #2 \rvert #3 \rvert #4 \rangle_{#5}}

\newcommand*\pfqskip{8mu} \catcode`,\active \newcommand*
\pfq{\begingroup\catcode`\,\active\def,{\mskip\pfqskip\relax}\dopfq}
\catcode`\,12
\def\dopfq#1#2#3#4#5{%
{}_{#1}\phi_{#2}\left(\genfrac{}{}{0pt}{}{#3}{#4}\,;\,#5\right)%
\endgroup}

 \catcode`,\active \newcommand*
\pFq{\begingroup\catcode`\,\active\def,{\mskip\pfqskip\relax}\dopFq}
\catcode`\,12
\def\dopFq#1#2#3#4#5{%
{}_{#1}F_{#2}\left(\genfrac{}{}{0pt}{}{#3}{#4}\,;\,#5\right)%
\endgroup}

\title{An algebraic interpretation \\of the multivariate $q$-Krawtchouk polynomials}
\author[$\dag$]{Vincent X. Genest}
\author[$\ddag$]{Sarah Post}
\author[*]{Luc Vinet}
\affil[$\dag$]{Department of Mathematics, Massachusetts Institute of Technology, 77 Massachusetts Avenue, Cambridge, MA 02139, USA}
\affil[$\ddag$]{Department of Mathematics, University of Hawai'i at Manoa, 2665 McCarthy Mall, Honolulu, HI 96822, USA}
\affil[*]{Centre de recherches math\'ematiques, Universit\'e de Montr\'eal, C.P.  6128 succ. Centre-ville, Montr\'eal (Qu\'ebec) H3C 3J7, Canada}
\date{}

\begin{document}
\maketitle
\thispagestyle{empty}
\hrule
\begin{abstract}\noindent
The multivariate quantum $q$-Krawtchouk polynomials are shown to arise as matrix elements of ``$q$-rotations'' acting on the state vectors of many $q$-oscillators. The focus is put on the two-variable case. The algebraic interpretation is used to derive the main properties of the poly\-nomials: orthogonality, duality, structure relations, difference equations and recurrence relations. The extension to an arbitrary number of variables is presented.\bigskip

\noindent\textbf{Keywords:} Multivariate $q$-Krawtchouk polynomials; $q$-oscillator algebra; $q$-rotations\smallskip

\noindent\textbf{AMS classification numbers:} 33D45, 16T05
\end{abstract}
\hrule
\section{Introduction}
The purpose of this paper is to provide an algebraic model for the multi-variable $q$-Kraw\-tchouk polynomials and to show how this model provides a cogent framework for the characterization of these ortho\-gonal functions. The algebraic interpretation presented here is in terms of matrix elements of unitary ``$q$-rotations'' acting on $q$-oscillator states; these $q$-rotations are expressed as $q$-exponentials in the generators of several independent $q$-oscillator algebras. For illustration purposes, the focus will be put on the two-variable case. The algebraic model will lead to a natural and explicit derivation of the main properties satisfied by the bivariate quantum $q$-Krawtchouk polynomials: orthogonality relation, duality property, structure relations, $q$-difference equations and recurrence relations. How this approach generalizes directly to an arbitrary number of variables will also be explained. 

The standard (univariate) Krawtchouk polynomials form one of the simplest families of hypergeometric orthogonal polynomials of the Askey scheme \cite[Ch. 9]{2010_Koekoek_&Lesky&Swarttouw}. These polynomials of degree $n$ in $x$, denoted by $\kappa_{n}(x;p,N)$, have the expression
\begin{align*}
\kappa_{n}(x;p,N)=(-1)^{n}(-N)_{n}\;\pFq{2}{1}{-n,-x}{-N}{\frac{1}{p}},\qquad n=0,1,\ldots,N,
\end{align*}
where $0<p<1$ is a parameter and $N$ is a positive integer\footnote{The polynomials $\kappa_{n}(x;p,N)$ agree with $K_n(x;p,N)$ in \cite[Section 9.11]{2010_Koekoek_&Lesky&Swarttouw} up to a normalization factor. In this paper,  the uppercase $K$ are reserved for the multivariate polynomials.}. In the above, ${}_pF_{q}$ stands for the generalized hypergeometric series and
\begin{align*}
(a)_n=a(a+1)\cdots (a+n-1),\qquad (a)_0:=1,
\end{align*}
is the shifted factorial, or Pochhammer symbol \cite[(1.1.2)]{2001_Andrews&Askey&Roy}. The polynomials $\kappa_{n}(x;p,N)$ are orthogonal with respect to the binomial distribution, i.e.
\begin{align*}
\sum_{x=0}^{N}\binom{N}{x}p^{x}(1-p)^{N-x}\kappa_{n}(x;p,N)\,\kappa_{m}(x;p,N)=\eta_{n}\delta_{nm},
\end{align*}
where $\eta_{n}\neq 0$ are normalization coefficients. The Krawtchouk polynomials have a well-known algebraic interpretation  as matrix elements of unitary irreducible representations of $\mathfrak{so}(3)$; in point of fact, most of their properties follow from that interpretation (see for example \cite{1982_Koornwinder_SIAMJMathAnal_13_1011}\cite[Section 6.8]{ 1991_Vilenkin&Klimyk}). 

In \cite{1991_Tratnik_JMathPhys_32_2065, 1991_Tratnik_JMathPhys_32_2337} Tratnik presented a multi-variable extension of the Askey tableau, giving for each family of the scheme explicit expressions of the polynomials as well as of their ortho\-gonality measures; special cases of these multivariate polynomials occurred earlier, e.g. \cite{Karlin-McGregor-1975}. Much later in \cite{2010_Geronimo&Illiev_ConstrApprox_31_417}, Geronimo and Iliev showed that Tratnik's families of orthogonal polynomials are all bispectral: they exhibited the recurrence relations and the eigenvalue equations that these polynomials satisfy. Despite the fact that the multivariate polynomials defined by Tratnik are expressed in terms of one-variable polynomials, they are non-trivial and arise in a number of different contexts, see for example \cite{2014_Genest&Vinet_JPhysA_47_455201, 2011_Griffiths&Spano_Bernoulli_1095, Karlin-McGregor-1975, 2015_Post_SIGMA_11_57,  1998_Rosengren_SIAMJMathAnal_30_233, 2007_Scarabotti_MethApplAnal_14_355}.  In two variables, the Krawtchouk polynomials defined by Tratnik have the expression
\begin{align*}
K_{n_1,n_2}(x_1,x_2;p_1,p_2;N)=\kappa_{n_1}(x_1;p_1, x_1+x_2)\,\kappa_{n_2}(x_1+x_2-n_1; p_2, N-n_1),\quad n_1,n_2=0,1,\ldots N,
\end{align*}
with $n_1+n_2\leq N$ and are orthogonal with respect to the trinomial distribution
\begin{align*}
\sum_{\substack{x_1,x_2=0\\ x_1+x_2\leq N}}^{N}\binom{N}{x_1,x_2}\left(\frac{p_1p_2}{1-p_2}\right)^{x_1}&\left(\frac{(1-p_1)p_2}{1-p_2}\right)^{x_2}(1-p_2)^{N} 
\\
&\times K_{n_1,n_2}(x_1,x_2)\,K_{m_1,m_2}(x_1,x_2)
=\eta_{n_1,n_2}\delta_{n_1m_1}\delta_{n_2m_2},
\end{align*}
where $\binom{N}{x_1,x_2}$ are the trinomial coefficients and where $p_1$, $p_2$ are parameters. In \cite{2013_Genest&Vinet&Zhedanov_JPhysA_46_505203}, it was shown that Tratnik's multivariate Krawtchouk polynomials arise in the matrix elements of the reducible unitary rotation group representations on oscillator states \footnote{Let us note that in \cite{2013_Genest&Vinet&Zhedanov_JPhysA_46_505203} and \cite{2010_Geronimo&Illiev_ConstrApprox_31_417}, a different definition of Tratnik's two-variable Krawtchouk polynomials was used. The aforementioned algebraic interpretation is valid for both definitions.}. It was also shown that Tratnik's Krawtchouk polynomials are in fact a special case of the multivariate Krawtchouk polynomials introduced by Griffiths in \cite{1971_Griffiths_AusJStat_13_27} and by Hoare and Rahman in \cite{2008_Hoare&Rahman_SIGMA_4_89}. Similarly to the one-variable case, the properties of these multivariate functions follow from their group-theoretical interpretation \cite{2013_Genest&Vinet&Zhedanov_JPhysA_46_505203}; see also \cite{2012_Iliev&Terwilliger_TransAmerMathsoc_364_4225}.

Basic analogs of Tratnik's multivariate orthogonal polynomials were introduced by Gasper and Rahman in \cite{2007_Gasper&Rahman_Ramanujan_13_389}. Like Tratnik, Gasper and Rahman explicitly defined their families of multivariate $q$-deformed polynomials as convoluted products of one-variable polynomials of the basic Askey scheme and proved their orthogonality relations. In \cite{2011_Iliev_TransAmerMathSoc_363_1577}, Iliev showed that these polynomials, which form a multivariate extension of the basic Askey tableau, are also bispectral. We note here that some of these families of orthogonal polynomials, namely the multivariate $q$-Hahn and $q$-Racah polynomials, have already been seen to arise in algebraic and combinatorial contexts \cite{2001_Rosengren_IJMMS_28_331, 2011_Scarabotti_RamanujanJ_25_57}. The two-variable $q$-Krawtchouk polynomials introduced in \cite{2007_Gasper&Rahman_Ramanujan_13_389} have the expression
\begin{align}
\label{2var-qKrawtchouk}
\mathbf{K}_{n_1,n_2}(x_1,x_2;\alpha_1,\alpha_2;N)=k_{n_1}(x_1;\alpha_1^{-2},x_1+x_2;q)\;k_{n_2}(x_1+x_2-n_1; \alpha_2^{-2}, N-n_1;q),
\end{align}
with $n_1,n_2$ are non-negative integers such that $n_1+n_2\leq N$. Here $k_{n}(x,p,N;q)$ stands for the one-variable quantum $q$-Krawtchouk polynomials of degree $n$ in $q^{-x}$ which are defined as \cite[Ch. 14]{2010_Koekoek_&Lesky&Swarttouw}
\begin{align}
\label{1var-qKraw}
k_{n}(x; p,N;q)=(-1)^{n}(q^{-N};q)_{n}q^{\binom{n}{2}}\;\pfq{2}{1}{q^{-n},q^{-x}}{q^{-N}}{q,\,p q^{n+1}},\quad n=0,1,\ldots, N,
\end{align}
where $0<q<1$, $p>q^{-N}$ and where $N$ is a positive integer\footnote{The polynomials $k_{n}(x; p,N;q)$ agree with $K_n^{qtm}(q^x; p, N, q)$ of \cite[Section 14.14]{2010_Koekoek_&Lesky&Swarttouw} up to a normalization factor. }. In the above
\begin{align*}
(a;q)_{n}=(1-a)(1-a q)\cdots (1-a q^{n-1}),\qquad (a;q)_0:=1,
\end{align*}
stands for the $q$-Pochhammer symbol and ${}_r\phi_{s}$ is the basic hypergeometric series defined as \cite[(1.2.22)]{2004_Gasper&Rahman}
\begin{align*}
{}_r\phi_{s}\left(\genfrac{}{}{0pt}{}{a_1, a_2, \ldots, a_{r}}{b_1,b_2,\ldots, b_{s}};\; q,\,z\right)=
 \sum_{k\geq 0}\frac{(a_1;q)_{k}(a_2;q)_{k}\cdots (a_{r};q)_{k}}{(q;q)_{k}(b_1;q)_{k}(b_2;q)_{k}\cdots (b_{s};q)_{k}}\left[(-1)^{k}q^{\binom{k}{2}}\right]^{1+s-r}z^{k}.
\end{align*}
In this paper, we shall provide an algebraic interpretation of Gasper and Rahman's multivariate (quantum) $q$-Krawtchouk polynomials and show how the properties of these polynomials can be obtained from that interpretation. On the one hand, this work can be considered as a multi-variable generalization of \cite{2015_Genest&Post&Vinet&Yu&Zhedanov}, where the algebraic interpretation of the one-variable quantum $q$-Krawtchouk in terms of $q$-rotations was investigated. On the other hand, it can be viewed as a $q$-generalization of \cite{2013_Genest&Vinet&Zhedanov_JPhysA_46_505203} (in the Tratnik special case), where the interpretation of the multivariate Krawtchouk polynomials in terms of rotation group representations on oscillator states was studied. For simplicity the emphasis shall be put on the two-variable case, with the understanding that the approach directly extends to an arbitrary number of variables.

The outline of the paper is as follows. In Section 2, the interpretation of the one-variable quantum $q$-Krawtchouk polynomials in terms of $q$-rotations, investigated in \cite{2015_Genest&Post&Vinet&Yu&Zhedanov}, is briefly reviewed. In Section 3, the algebraic model for the two-variable quantum $q$-Krawtchouk polynomials is presented. The bivariate quantum $q$-Krawtchouk are seen to arise as matrix elements of unitary three-dimensional $q$-rotations acting on the state space of three mutually commuting copies of the $q$-oscillator algebra. The orthogonality weight and the duality relation for the polynomials are derived. In Section 4, new structure relations for two-variable $q$-Krawtchouk polynomials are found using the algebraic setting. In Section 5, the recurrence relations and the difference equations are obtained and are seen to coincide with those found by Iliev in \cite{2011_Iliev_TransAmerMathSoc_363_1577}. The multivariate case is discussed in the conclusion.
\section{Review of the univariate case}
In this section, the algebraic interpretation of the one-variable quantum $q$-Krawtchouk polynomials as matrix elements of unitary $q$-rotations acting on the space of two $q$-oscillators is reviewed. Mild generalizations of the $q$-Baker--Campbell--Hausdorff formulas are also presented.
\subsection{The $q$-oscillator algebra and the Schwinger realization of $U_{q}(\mathfrak{sl}_2)$}
Let $0<q<1$ and consider two mutually commuting copies of the $q$-oscillator algebra with generators $A_{\pm}$, $A_0$ and $B_0$, $B_{\pm}$ that satisfy the relations
\begin{alignat}{3}
\label{qOsc-Algebra}
\begin{aligned}
[X_0, X_{\pm}]&=\pm X_{\pm},\qquad& [X_{-},X_{+}]&=q^{X_0},\qquad& X_{-}X_{+}-q X_{+}X_{-}&=1,
\end{aligned}
\end{alignat}
with $X\in\{A,B\}$ and where $[X\_,Y\_]=0$ if $X\neq Y$. The algebra \eqref{qOsc-Algebra} has a standard infinite-dimensional representation on the ortho\-normal basis states
\begin{align*}
\ket{n_{A},n_{B}}:=\ket{n_{A}}\otimes \ket{n_{B}},
\end{align*}
where $n_{A}$ and $n_{B}$ are non-negative integers. This representation is defined by the following action of the generators on the factors of the direct product states:
\begin{align}
\label{Action}
X_{+}\ket{n_{X}}&=\sqrt{\frac{1-q^{n_{X}+1}}{1-q}} \ket{n_{X}+1},\qquad X_{-}\ket{n_{X}}=\sqrt{\frac{1-q^{n_{X}}}{1-q}}\ket{n_{X}-1},\qquad 
X_0\ket{n_{X}}=n_{X}\ket{n_{X}},
\end{align}
where $X=A \text{ or } B$. Note that in this representation, $X_{+}$ and $X_{-}$ are mutual adjoints and $X_0$ is self-adjoint. It is easily seen that in the $q\uparrow 1$ limit, the representation defined by the actions \eqref{Action} goes to the standard oscillator representation.

The Schwinger realization of the quantum algebra $U_{q}(\mathfrak{sl}_2)$ is constructed from the two $q$-oscillator algebras \eqref{qOsc-Algebra} by taking  \cite{1989_Biedenharn_JPhysA_22_L873}
\begin{align*}
J_{+}=q^{-\frac{A_0+B_0-1}{4}} A_{+}B_{-},\qquad J_{-}=q^{-\frac{A_0+B_0-1}{4}}A_{-}B_{+},\qquad J_0=\frac{A_0-B_0}{2}.
\end{align*}
Indeed, it is easily verified using the relations \eqref{qOsc-Algebra} that the operators $J_{\pm}$ and $J_0$ satisfy the defining relations of $U_{q}(\mathfrak{sl}_2)$ which read
\begin{align*}
[J_0, J_{\pm}]=\pm J_{\pm},\qquad [J_{+},J_{-}]=\frac{q^{J_0}-q^{-J_0}}{q^{1/2}-q^{-1/2}}.
\end{align*}
As can be verified directly using \eqref{Action}, the orthonormal $q$-oscillator states defined by
\begin{align}
\label{Basis-1}
\kket{n}{N}:=\ket{n,N-n},\qquad n=0,1,\ldots, N,
\end{align}
where $N$ is a non-negative integer, support the standard unitary $(N+1)$-dimensional irreducible representations of $U_{q}(\mathfrak{sl}_2)$.
\subsection{$q$-analogs of the BCH formulas and $q$-rotation operators}
Consider the little and big $q$-exponential functions, respectively denoted by $e_{q}(z)$ and $E_{q}(z)$, which are defined as follows:
\begin{align}
\label{q-Exp}
e_{q}(z)=\sum_{n=0}^{\infty}\frac{z^{n}}{(q;q)_{n}}=\frac{1}{(z;q)_{\infty}},\quad |z|<1,\qquad \text{and}\qquad E_{q}(z)=\sum_{n=0}^{\infty}\frac{q^{\binom{n}{2}}z^{n}}{(q;q)_{n}}=(-z;q)_{\infty}.
\end{align}
One clearly has $e_{q}(z)E_{q}(-z)=1$. In the following, two $q$-analogs of the Baker--Campbell--Hausdorff relation shall be needed.  The first relation is of the form
\begin{align}
\label{BCH-1}
E_{q}(\lambda X) Y e_{q}(-\lambda\,q^{\alpha} X)=\sum_{n=0}^{\infty}\frac{\lambda^{n}}{(q;q)_{n}}[X,Y]_{n},
\end{align}
where $[X,Y]_{n}$ is defined recursively through
\begin{align*}
[X,Y]_0=1,\qquad [X,Y]_{n+1}=q^{n} X\,[X,Y]_{n}-q^{\alpha}\,[X,Y]_{n}\,X,\qquad n=1,2,\ldots.
\end{align*}
The second relation reads
\begin{align}
\label{BCH-2}
e_{q}(\lambda X) Y E_{q}(-\lambda\,q^{\alpha} X)=\sum_{n=0}^{\infty} \frac{\lambda^{n}}{(q;q)_{n}}[X,Y]_{n}',
\end{align}
where $[X,Y]_{n}'$  is defined recursively by
\begin{align*}
[X,Y]_0'=1,\qquad [X,Y]_{n+1}'=X\,[X,Y]_{n}'-q^{n+\alpha}\,[X,Y]_{n}'\,X,\qquad n=0,1,2,\ldots.
\end{align*}
The identities \eqref{BCH-1}, \eqref{BCH-2} can be verified in a straightforward manner by expanding the $q$-exponentials in power series using \eqref{q-Exp}. For $\alpha=0$, these relations are given in \cite{1993_Floreanini&Vinet_PhysLettA_180_393}, among others.

Let $\theta$ be a real number such that $|\theta|<1$ and let $U_{AB}(\theta)$ be the operator defined as
\begin{multline}
\label{UAB}
U_{AB}(\theta)=
\\
e_{q}^{1/2}\left(\theta^2 q^{-A_0}\right)\,e_{q}\left(\theta(1-q)q^{-(A_0+B_0)/2}A_{+}B_{-}\right)E_{q}\left(-\theta(1-q)q^{-(A_0+B_0)/2}A_{-}B_{+}\right)\,E_{q}^{1/2}\left(-\theta^2 q^{-B_0}\right).
\end{multline}
This operator will be called a ``$q$-rotation'' operator. This term, coined by Zhedanov in \cite{1993_Zhedanov_JMathPhys_34_2631}, comes from the fact that $U_{AB}(\theta)$ can be viewed as a $q$-analog of an $SU(2)$ element obtained via the exponential map from the algebra to the group; see \cite{2015_Genest&Post&Vinet&Yu&Zhedanov, 1993_Zhedanov_JMathPhys_34_2631} for more details. As shown in \cite{2015_Genest&Post&Vinet&Yu&Zhedanov}, the operator \eqref{UAB} is unitary: it satisfies the relations $U^{\dagger}U=1$ and $UU^{\dagger}=1$. The unitarity property of $U_{AB}(\theta)$ follows directly from the formulas
\begin{align*}
e_{q}(\alpha\,A_{-}B_{+})\,e_{q}\left(\frac{\alpha\beta}{(1-q)^2}q^{B_0}\right)\,e_{q}(\beta A_{+}B_{-})=e_{q}(\beta A_{+}B_{-})e_{q}\left(\frac{\alpha\beta}{(1-q)^2}q^{A_0}\right)\,e_{q}(\alpha A_{-}B_{+}),
\end{align*}
and 
\begin{align*}
E_{q}(\gamma A_{+}B_{-})E_{q}\left(-\frac{\gamma\delta}{(1-q)^2}q^{B_0}\right)E_{q}(\delta A_{-}B_{+})=E_{q}(\delta A_{-}B_{+})E_{q}\left(-\frac{\gamma \delta}{(1-q)^2}q^{A_0}\right)E_{q}(\gamma A_{+}B_{-}),
\end{align*}
where $\alpha,\beta,\gamma,\delta$ are constants or central elements. These formulas can be proved using the $q$-BCH relations \eqref{BCH-1} and \eqref{BCH-2} as well as elementary properties of the $q$-exponential functions; see \cite{2015_Genest&Post&Vinet&Yu&Zhedanov}.
\subsection{Matrix elements and univariate quantum $q$-Krawtchouk polynomials}
Let us now recall how the one-variable quantum $q$-Krawtchouk polynomials arise in this setting. The matrix elements of the unitary $q$-rotation operator \eqref{UAB} in the basis \eqref{Basis-1} are defined as
\begin{align}
\label{Univariate-Elements}
\xi_{n,x}^{(N)}(\theta)=\BBraket{N}{n}{U_{AB}(\theta)}{x}{N},
\end{align}
where $n$ and $x$ take values in $\{0,1,\ldots, N\}$. The matrix elements \eqref{Univariate-Elements} can be written as \cite{2015_Genest&Post&Vinet&Yu&Zhedanov}
\begin{align}
\label{1var-Explicit}
\xi_{n,x}(\theta)=\omega_{x}^{(N)}(\theta)\,\sigma_{n}^{(N)}(\theta)\,k_{n}(x,\theta^{-2},N;q),
\end{align}
where $k_{n}(x,p,N;q)$ are the quantum $q$-Krawtchouk polynomials defined in \eqref{1var-qKraw}. In \eqref{1var-Explicit}, the coefficient $\sigma_{n}^{(N)}$ has the expression
\begin{align}
\sigma_{n}^{(N)}(\theta)=\frac{(-1)^{n} q^{-\binom{n}{2}}}{(q^{-N};q)_{n}}\sqrt{\binom{N}{n}_{q}\frac{\theta^{2n}q^{-nN}}{(\theta^2 q^{-n};q)_{n}}},
\end{align}
where 
\begin{align*}
\binom{N}{n}_{q}=\frac{(q;q)_{N}}{(q;q)_{n}(q;q)_{N-n}},
\end{align*}
stands for the $q$-binomial coefficient. The expression for $\omega_{x}^{(N)}(\theta)$, which is defined as the matrix element
$\BBraket{N}{0}{U_{AB}(\theta)}{x}{N}$, reads
\begin{align}
\omega_{x}^{(N)}(\theta):=\BBraket{N}{0}{U_{AB}(\theta)}{x}{N}=(-1)^{x}q^{\binom{x}{2}}\sqrt{\binom{N}{x}_{q}\frac{(\theta^2 q^{-N};q)_{N}}{(\theta^2 q^{-N};q)_{x}}\theta^{2x}q^{-xN}}.
\end{align}
From the algebraic interpretation \eqref{1var-Explicit} of the one-variable quantum $q$-Krawtchouk polynomials, one can derive the main properties of these polynomials  \cite{2015_Genest&Post&Vinet&Yu&Zhedanov}. In what follows, the duality relation
\begin{align}
\label{duality}
\xi_{n,x}^{(N)}(\theta)=\xi_{N-x,N-n}^{(N)}(\theta),
\end{align}
which follows from the reality of the matrix elements \eqref{Univariate-Elements}, the identity $U^{\dagger}_{AB}(\theta)=U^{-1}_{AB}(\theta)$ and the observation that $U^{-1}_{AB}(\theta)=U_{BA}(\theta)$, shall prove particularly useful.
\section{A model for the two-variable $q$-Krawtchouk polynomials}
In this section, the algebraic model for the two-variable $q$-Krawtchouk polynomials, which involves a three-dimensional $q$-rotation acting on the states of three $q$-oscillators, is constructed. The orthogonality relation and the duality property of the polynomials are derived from the model.
\subsection{The model}
Consider the algebra generated by three mutually commuting $q$-oscillators with generators $X_0$, $X_{\pm}$, $X\in\{A,B,C\}$, satisfying the commutation relations \eqref{qOsc-Algebra} and consider its representation defined by the actions \eqref{Action} on the three-fold tensor product space spanned by the orthonormal basis vectors $\ket{n_A, n_{B}, n_{C}}=\ket{n_{A}}\otimes \ket{n_{B}}\otimes \ket{n_{C}}$, where $n_{A}$, $n_{B}$ and $n_{C}$ are non-negative integers. Let $N$ be a positive integer and introduce the basis states $\kket{n_1,n_2}{N}$ defined as
\begin{align}
\label{2var-Basis}
\kket{n_1,n_2}{N}:=\ket{n_1,n_2,N-n_1-n_2},\qquad n_1,n_2\in \{0,1,\ldots, N\},\qquad n_1+n_2\leq N.
\end{align}
For a given $N$, the states \eqref{2var-Basis} support an $(N+1)(N+2)/2$-dimensional irreducible representation of the quantum algebra $U_{q}(\mathfrak{sl}_3)$ realized with three independent $q$-oscillators; see \cite{1990_Hayashi_CommMathPhys_127_129}. Following \eqref{UAB}, the unitary operator $U_{XY}(\theta)$ effecting the $q$-rotation in the $XY$ ``plane'' is defined as 
\begin{multline}
\label{UXY}
U_{XY}(\theta)=
\\
e_{q}^{1/2}\left(\theta^2 q^{-X_0}\right)\,e_{q}\left(\theta(1-q)q^{-(X_0+Y_0)/2}X_{+}Y_{-}\right)E_{q}\left(-\theta(1-q)q^{-(X_0+Y_0)/2}X_{-}Y_{+}\right)\,E_{q}^{1/2}\left(-\theta^2 q^{-Y_0}\right).
\end{multline}
\subsection{Matrix elements and bivariate $q$-Krawtchouk polynomials}
Consider the matrix elements
\begin{align}
\label{2var-MatrixElements}
\Xi_{n_1,n_2}^{(N)}(x_1,x_2; \theta,\phi)=\BBraket{N}{n_1,n_2}{U_{BC}(\phi)\;U_{AB}(\theta)}{x_1,x_2}{N},
\end{align}
of the operator associated to a particular $q$-rotation in three dimensions. To ease the notation, $\Xi_{n_1,n_2}^{(N)}(x_1,x_2; \theta,\phi)$ shall sometimes be written simply as $\Xi_{n_1,n_2}^{(N)}(x_1,x_2)$. An explicit expression for the matrix elements \eqref{2var-MatrixElements} can be obtained by inserting a resolution of the identity in between the two $q$-rotation operators $U_{BC}(\phi)$ and $U_{AB}(\theta)$. This leads to
\begin{align}
\label{Decompo}
\Xi_{n_1,n_2}^{(N)}(x_1,x_2; \theta,\phi)=\sum_{\substack{s,t=0\\ s+t\leq N}}^{N}
\BBraket{N}{n_1,n_2}{U_{BC}(\phi)}{s,t}{N}\;\BBraket{N}{s,t}{U_{AB}(\theta)}{x_1,x_2}{N}.
\end{align}
Since the $q$-rotation operator $U_{BC}(\phi)$ leaves invariant the first index when acting on the basis vector $\kket{n_1,n_2}{N}$, it follows from the definition \eqref{Univariate-Elements} and the orthonormality of the basis that
\begin{align}
\label{First}
\BBraket{N}{n_1,n_2}{U_{BC}(\phi)}{s,t}{N}=\delta_{n_1, s}\,\xi_{n_2,t}^{(N-n_1)}(\phi).
\end{align}
Similarly, when acting on the state $\kket{x_1,x_2}{N}$, the operator $U_{AB}(\theta)$ leaves invariant the sum of the first two indices $x_1+x_2$ and one can thus write
\begin{align}
\label{Second}
\BBraket{N}{s,t}{U_{AB}(\theta)}{x_1,x_2}{N}=\delta_{s+t, x_1+x_2 }\;\xi_{s, x_1}^{(x_1+x_2)}(\theta).
\end{align}
Upon using \eqref{First} and \eqref{Second} in \eqref{Decompo}, the double summation cancels and one finds
\begin{align}
\label{Convo}
\Xi_{n_1,n_2}^{(N)}(x_1,x_2;\theta,\phi)=\xi_{n_1,x_1}^{(x_1+x_2)}(\theta)\;\xi_{n_2,x_1+x_2-n_1}^{(N-n_1)}(\phi).
\end{align}
With the explicit expression \eqref{1var-Explicit} of the one-dimensional $q$-rotation matrix elements in terms of the univariate quantum $q$-Krawtchouk polynomials, one has 
\begin{multline}
\label{Exp-1}
\Xi_{n_1,n_2}^{(N)}(x_1,x_2;\theta,\phi)=\omega_{x_1}^{(x_1+x_2)}(\theta)\,\sigma_{n_1}^{(x_1+x_2)}(\theta)\,k_{n_1}(x_1, \theta^{-2}, x_1+x_2;q)
\\
\times \omega_{x_1+x_2-n_1}^{(N-n_1)}(\phi)\,\sigma_{n_2}^{(N-n_1)}(\phi)\;k_{n_2}(x_1+x_2-n_1,\phi^{-2},N-n_1;q).
\end{multline}
Comparing \eqref{Exp-1} with \eqref{2var-qKrawtchouk}, it is seen that the matrix elements $\Xi_{n_1,n_2}^{(N)}(x_1,x_2)$ can be expressed in terms of the two-variable $q$-Krawtchouk polynomials $\mathbf{K}_{n_1,n_2}(x_1,x_2;\alpha_1,\alpha_2;N)$ with parameters $\alpha_1=\theta$ and $\alpha_2=\phi$. One must however verify that the matrix elements $\Xi_{n_1,n_2}^{(N)}(x_1,x_2)$ factorize as a product of the matrix element $\BBraket{N}{0,0}{U_{BC}(\phi)U_{AB}(\theta)}{x_1,x_2}{N}$, a normalization factor depending only on the degree indices $n_1$ and $n_2$, and the polynomials $\mathbf{K}_{n_1,n_2}(x_1,x_2;\alpha_1,\alpha_2;N)$ themselves. A direct, though technical, calculation shows that this is indeed the case and that one can write the matrix elements $\Xi_{n_1,n_2}^{(N)}(x_1,x_2)$ as
\begin{align}
\label{2var-Elements-Explicit}
\Xi_{n_1,n_2}^{(N)}(x_1,x_2;\theta,\phi)=W_{x_1,x_2}^{(N)}(\theta,\phi)\,\Sigma_{n_1,n_2}^{(N)}(\theta,\phi)\,\mathbf{K}_{n_1,n_2}(x_1,x_2;\theta,\phi;N),
\end{align}
where $W_{x_1,x_2}^{(N)}(\theta,\phi)$ is given by
\begin{multline}
\label{W}
W_{x_1,x_2}^{(N)}(\theta,\phi)=\BBraket{N}{0,0}{U_{BC}(\phi)U_{AB}(\theta)}{x_1,x_2}{N}
\\
=(-1)^{x_2}\Bigg[(-1)^{N-x_1}\,(\theta^{2}q^{-1})^{x_1+x_2}\,q^{\binom{x_1}{2}}\;\binom{N}{x_1,x_2}_{q}
 (\theta^{-2}q;q)_{x_2}(\phi^{-2}q;q)_{N-x_1-x_2}\;\phi^{2N}q^{-N(N+1)/2}
 \Bigg]^{1/2},
\end{multline}
where the normalization factor $\Sigma_{n_1,n_2}^{(N)}(\theta,\phi)$ has the expression
\begin{multline}
\label{Sigma}
\Sigma_{n_1,n_2}^{(N)}(\theta,\phi)=q^{-\binom{n_1}{2}-\binom{n_2}{2}}\,(-1)^{n_1}\,\phi^{-n_1}\frac{(q;q)_{N-n_1-n_2}}{(q;q)_{N}}
\\
\times
 \Bigg[
 (-1)^{n_1+n_2}
 \binom{N}{n_1,n_2}_{q}\frac{q^{N(n_1+n_2)+2n_1+n_2-n_1n_2-\binom{n_1}{2}-\binom{n_2}{2}}}{(q\phi^{-2};q)_{n_2}(q\theta^{-2};q)_{n_1}}
 \Bigg]^{1/2},
\end{multline}
and where $\mathbf{K}_{n_1,n_2}(x_1,x_2;\alpha_1,\alpha_2;N)$ are the two-variable quantum $q$-Krawtchouk polynomials defined in \eqref{2var-qKrawtchouk}. In \eqref{W} and \eqref{Sigma}, $\binom{N}{x,y}_{q}$ stands for the $q$-trinomial coefficients
\begin{align*}
\binom{N}{x,y}_{q}=\frac{(q;q)_{N}}{(q;q)_{x}(q;q)_{y}(q;q)_{N-x-y}}.
\end{align*}
\textbf{Remark.} Let us note here that in \eqref{2var-MatrixElements}, the choice of the ordering of the two $q$-rotations to form a $q$-rotation in three space is important. For example, consider the matrix elements 
\begin{align*}
\BBraket{N}{n_1,n_2}{U_{AC}(\phi)U_{BC}(\theta)}{x_1,x_2}{N}=\xi_{n_1,x_1}^{(N-n_2)}(\theta)\,\xi_{n_2,x_2}^{(N-x_1)}(\phi),
\end{align*}
which can be seen as $q$-analogs of the matrix elements leading to Geronimo and Iliev's two-variable Krawtchouk polynomials \cite{2010_Geronimo&Illiev_ConstrApprox_31_417}. It can be verified directly that with this choice of ordering one has a factorization the form
\begin{align*}
\BBraket{N}{n_1,n_2}{U_{AC}(\phi)U_{BC}(\theta)}{x_1,x_2}{N}=\BBraket{N}{0,0}{U_{AC}(\phi)U_{BC}(\theta)}{x_1,x_2}{N}\;P_{n_1,n_2}(x_1,x_2),
\end{align*}
where $P_{n_1,n_2}(x_1,x_2)$, with $P_{0,0}(x_1,x_2)=1$, are not polynomial in the proper variables.
\subsection{Orthogonality relations}
The orthogonality relation for the two-variable $q$-Krawtchouk polynomials immediately follows from the unitarity of the $q$-rotation $U_{BC}(\phi)U_{AB}(\theta)$, the orthonormality of the basis states and the reality of the matrix elements \eqref{2var-MatrixElements}. One has
\begin{align}
\label{Ortho-1}
\begin{aligned}
\bbraket{N}{n_1,n_2}{m_1,m_2}{N}&=\BBraket{N}{n_1,n_2}{U_{BC}(\phi)U_{AB}(\theta)U_{AB}^{\dagger}(\theta)U_{BC}^{\dagger}(\phi)}{m_1,m_2}{N}
\\
 &=\sum_{\substack{x_1,x_2=0\\ x_1+x_2\leq N}}^{N}
\BBraket{N}{n_1,n_2}{U_{BC}(\phi)U_{AB}(\theta)}{x_1,x_2}{N}\;\BBraket{N}{x_1,x_2}{U_{AB}^{\dagger}(\theta)U_{BC}^{\dagger}(\phi)}{m_1,m_2}{N}
\\
&=\sum_{\substack{x_1,x_2=0\\ x_1+x_2\leq N}}^{N} \Xi_{n_1,n_2}^{(N)}(x_1,x_2;\theta,\phi)\;\Xi_{m_1,m_2}^{(N)}(x_1,x_2;\theta,\phi)=\delta_{n_1 m_1}\delta_{n_2m_2}.
\end{aligned}
\end{align}
Upon introducing \eqref{2var-Elements-Explicit} in the above relation, one finds that the two-variable quantum $q$-Krawtchouk polynomials satisfy the orthogonality relation
\begin{align}
\label{Ortho}
\sum_{\substack{x_1,x_2=0\\ x_1+x_2\leq N}}^{N} \Delta_{x_1,x_2}^{(N)}(\theta,\phi)\;\mathbf{K}_{n_1,n_2}(x_1,x_2;\theta,\phi;N)\mathbf{K}_{m_1,m_2}(x_1,x_2;\theta,\phi;N)
=H_{n_1,n_2}^{(N)}(\theta,\phi)\delta_{n_1 m_1}\delta_{n_2 m_2},
\end{align}
where $\Delta(x_1,x_2)^{(N)}(\theta,\phi)=[W_{x_1,x_2}^{(N)}(\theta,\phi)]^2$ and $H_{n_1,n_2}^{(N)}(\theta,\phi)=[\Sigma_{n_1,n_2}^{(N)}(\theta,\phi)]^{-2}$. From \eqref{Ortho-1}, it is easily seen that one can also write a dual orthogonality relation of the form
\begin{align}
\label{Ortho-Dual}
\sum_{\substack{n_1,n_2=0\\ n_1+n_2\leq N}}^{N} \widetilde{\Delta}_{n_1,n_2}^{(N)}(\theta,\phi)\mathbf{K}_{n_1,n_2}(x_1,x_2;\theta,\phi;N)\mathbf{K}_{n_1,n_2}(y_1,y_2;\theta,\phi;N)
=\widetilde{H}_{x_1,x_2}^{(N)}(\theta,\phi)\delta_{x_1 y_1}\delta_{x_2 y_2},
\end{align}
where $\widetilde{\Delta}_{n_1,n_2}^{(N)}(\theta,\phi)=[\Sigma_{n_1,n_2}^{(N)}(\theta,\phi)]^{2}$ and $\widetilde{H}_{x_1,x_2}^{(N)}=[W_{x_1,x_2}^{(N)}(\theta,\phi)]^{-2}$.
\subsection{Duality}
The matrix elements \eqref{2var-MatrixElements} satisfy a duality relation that allows to interchange the roles of the degree indices $n_1$ and $n_2$ with those of the variables indices $x_1$ and $x_2$. This duality relation is obtained by combining the duality relation \eqref{duality} and the expression \eqref{Convo} for the matrix elements. One has
\begin{align}
\label{duality-2}
\begin{aligned}
\Xi_{n_1,n_2}^{(N)}(x_1,x_2;\theta,\phi)&=\xi_{n_1,x_1}^{(x_1+x_2)}(\theta)\;\xi_{n_2,x_1+x_2-n_1}^{(N-n_1)}(\phi)
=\xi_{x_2,\;x_1+x_2-n_1}^{(x_1+x_2)}(\theta)\xi_{N-x_1-x_2,\;N-n_1-n_2}^{(N-n_1)}(\phi)
\\
&=\xi_{N-x_1-x_2,\;N-n_1-n_2}^{(N-n_1)}(\phi)\xi_{x_2,\;x_1+x_2-n_1}^{(x_1+x_2)}(\theta)=\Xi_{N-x_1-x_2,\,x_2}^{(N)}(N-n_1-n_2, n_2; \phi, \theta),
\end{aligned}
\end{align}
where we have used \eqref{Convo} in the first step and \eqref{duality} in the second step; note that the parameters $\theta$ and $\phi$ are interchanged in the duality transformation. In terms of the polynomials $\mathbf{K}_{n_1,n_2}(x_1,x_2;\theta,\phi;N)$, the duality relation \eqref{duality-2} becomes
\begin{multline}
(-1)^{x_2+n_2}q^{x_1^2+x_2^2+x_1x_2+x_1-N(x_1+x_2-1)} \left(\frac{(q\theta^{-2};q)_{x_2}(q \phi^{-2};q)_{N-x_1-x_2}}{(q;q)_{x_1}}\right) \mathbf{K}_{n_1,n_2}(x_1,x_2;\theta,\phi;N)
\\
=q^{n_1^2+n_2^2+n_1n_2-n_1-n_2-N(n_1+n_2)}\left(\frac{(q\theta^{-2};q)_{n_1}(q \phi^{-2};q)_{n_2}}{(q;q)_{N-n_1-n_2}}\right)\mathbf{K}_{N-x_1-x_2,x_2}(N-n_1-n_2, n_2;\phi,\theta;N).
\end{multline}

\section{Structure relations}
In this section, the $q$-BCH formulas \eqref{BCH-1}, \eqref{BCH-2} and  the algebraic model constructed in the previous section are used to obtain structure relations for the two-variable quantum $q$-Krawtchouk polynomials. Raising and lowering relations on the degrees and the variables are found.
\subsection{Structure relations with respect to the index $n_1$}
We first use the algebraic properties of the three-dimensional $q$-rotation operator $U_{BC}(\phi)U_{AB}(\theta)$ to obtain structures relations on the first degree index $n_1$.
\subsubsection{Raising relation in $n_1$}
To get a raising relation with respect to the first degree index $n_1$, consider the matrix element $\BBraket{N-1}{n_1,n_2}{A_{-} U_{BC}(\phi)U_{AB}(\theta)}{x_1,x_2}{N}$. On the one hand, using \eqref{Action} and \eqref{2var-MatrixElements}, one has
\begin{align}
\label{rr-1}
\BBraket{N-1}{n_1,n_2}{A_{-} U_{BC}(\phi)U_{AB}(\theta)}{x_1,x_2}{N}=\sqrt{\frac{1-q^{n_1+1}}{1-q}}\;\Xi_{n_1+1,n_2}^{(N)}(x_1,x_2;\theta,\phi).
\end{align}
On the other hand, one can write
\begin{align}
\label{rr-2}
\BBraket{N-1}{n_1,n_2}{A_{-} U_{BC}(\phi)U_{AB}(\theta)}{x_1,x_2}{N}
=\BBraket{N-1}{n_1,n_2}{U_{BC}(\phi)U_{AB}(q^{-1/2}\theta) \;\widetilde{A}_{-}}{x_1,x_2}{N},
\end{align}
where $\widetilde{A}_{-}$ is of the form
\begin{align}
\label{A-Tilde}
\widetilde{A}_{-}=U_{AB}^{\dagger}(q^{-1/2}\theta)U_{BC}^{\dagger}(\phi)\;A_{-}\;U_{BC}(\phi) U_{AB}(\theta)=U_{AB}^{\dagger}(q^{-1/2}\theta)\;A_{-}\;U_{AB}(\theta).
\end{align}
A direct calculation using the expression \eqref{UXY} and the formulas \eqref{BCH-1} and \eqref{BCH-2} yields
\begin{align}
\label{Rel1}
U_{AB}^{\dagger}(q^{-1/2}\theta)\; A_{-}\; U_{AB}(\theta)=\sqrt{1-\theta^2q^{-1}q^{-B_0}}\,A_{-}+\theta q^{-1/2}q^{(A_0-B_0)/2}B_{-}.
\end{align}
To obtain \eqref{Rel1} from the $q$-Baker--Campbell--Hausdorff relations, one should recall that $U_{XY}^{\dagger}=U^{-1}_{XY}$. Upon substituting \eqref{Rel1} in  \eqref{rr-2}, using the actions \eqref{Action} and comparing with \eqref{rr-1}, one finds 
\begin{multline}
\label{Raising-Degree-1}
\sqrt{1-q^{n_1+1}}\,\Xi_{n_1+1,n_2}^{(N)}(x_1,x_2;\theta,\phi)=\sqrt{1-\theta^2 q^{-x_2-1}}\sqrt{1-q^{x_1}}\,\Xi_{n_1,n_2}^{(N-1)}(x_1-1,x_2;q^{-1/2}\theta, \phi)
\\
+\theta q^{(x_1-x_2)/2}\sqrt{1-q^{x_2}} \;\Xi_{n_1,n_2}^{(N-1)}(x_1,x_2-1;q^{-1/2}\theta,\phi).
\end{multline}
Using the expression \eqref{2var-Elements-Explicit} of the matrix elements in terms of the two-variable $q$-Krawtchouk polynomials, the relation \eqref{Raising-Degree-1} takes the form
\begin{align}
\label{rr-22}
\begin{aligned}
q^{-n_1}\mathbf{K}_{n_1+1,n_2}(x_1,x_2;q^{1/2}\theta,\phi;N)&=q^{-(x_1+x_2)}(1-q^{x_1})(1-\theta^{-2}q^{x_2})\;\mathbf{K}_{n_1,n_2}(x_1-1,x_2;\theta,\phi;N-1)
\\
&\qquad -(1-q^{-x_2})\;\mathbf{K}_{n_1,n_2}(x_1,x_2-1;\theta,\phi;N-1).
\end{aligned}
\end{align}
\subsubsection{Lowering relation in $n_1$}
A lowering relation on the first degree index $n_1$ can be obtained by considering instead the matrix element $\BBraket{N+1}{n_1,n_2}{A_{+} U_{BC}(\phi)U_{AB}(q^{-1/2}\theta)}{x_1,x_2}{N}$. One has
\begin{align}
\label{rr-7}
\BBraket{N+1}{n_1,n_2}{A_{+} U_{BC}(\phi)U_{AB}(q^{-1/2}\theta)}{x_1,x_2}{N}=\sqrt{\frac{1-q^{n_1}}{1-q}}\;\Xi_{n_1-1,n_2}^{(N)}(x_1,x_2;q^{-1/2}\theta,\phi),
\end{align}
and also
\begin{align}
\label{rr-6}
\BBraket{N+1}{n_1,n_2}{A_{+} U_{BC}(\phi)U_{AB}(q^{-1/2}\theta)}{x_1,x_2}{N}=
\BBraket{N+1}{n_1,n_2}{U_{BC}(\phi)U_{AB}(\theta) \widetilde{A}_{+}}{x_1,x_2}{N},
\end{align}
with
\begin{align}
\label{Ap-Tilde}
\widetilde{A}_{+}=U_{AB}^{\dagger}(\theta)U_{BC}^{\dagger}(\phi) A_{+} U_{BC}(\phi) U_{AB}(q^{-1/2}\theta).
\end{align}
It is directly seen that \eqref{Ap-Tilde} is the complex conjugate of \eqref{A-Tilde}. Consequently, we have
\begin{align}
\label{Ap-Tilde-2}
\widetilde{A}_{+}=A_{+}\sqrt{1-\theta^2q^{-1}q^{-B_0}}+\theta q^{-1/2}B_{+}q^{(A_0-B_0)/2}.
\end{align}
Upon substituting \eqref{Ap-Tilde-2} in \eqref{rr-6}, using the actions \eqref{Action} and comparing with \eqref{rr-7}, one finds 
\begin{multline}
\sqrt{1-q^{n_1}}\Xi_{n_1-1,n_2}^{(N)}(x_1,x_2;\theta,\phi)=\sqrt{1-q^{x_1+1}}\sqrt{1-\theta^2q^{-x_2}}\;\Xi_{n_1,n_2}^{(N+1)}(x_1+1,x_2;q^{1/2}\theta,\phi)
\\
+\theta q^{(x_1-x_2)/2} \sqrt{1-q^{x_2+1}} \;\Xi_{n_1,n_2}^{(N+1)}(x_1,x_2+1;q^{1/2}\theta,\phi),
\end{multline}
where we have taken $\theta\rightarrow q^{1/2}\theta$. Using \eqref{2var-Elements-Explicit}, this gives 
\begin{multline}
\frac{q^{n_1}(1-q^{n_1})}{\theta^{2}}\;\mathbf{K}_{n_1-1,n_2}(x_1,x_2;\theta,\phi;N)=q^{x_1+2}\;\mathbf{K}_{n_1,n_2}(x_1,x_2+1;q^{1/2}\theta,\phi;N+1)
\\
-q^{x_1+2}\mathbf{K}_{n_1,n_2}(x_1+1,x_2;q^{1/2}\theta,\phi;N+1).
\end{multline}

\subsection{Structure relations with respect to the index $n_2$}
We now look for structure relations on the second degree index $n_2$.
\subsubsection{Raising relation in $n_2$}
To get a raising relation with respect the second degree index $n_2$, consider the matrix element $\BBraket{N-1}{n_1,n_2}{q^{A_0/2}B_{-}U_{BC}(\phi)U_{AB}(\theta)}{x_1,x_2}{N}$. One has
\begin{align}
\label{rr-3}
\BBraket{N-1}{n_1,n_2}{q^{A_0/2}B_{-}U_{BC}(\phi)U_{AB}(\theta)}{x_1,x_2}{N}=q^{n_1/2}\sqrt{\frac{1-q^{n_2+1}}{1-q}}\,\Xi_{n_1,n_2+1}^{(N)}(x_1,x_2;\theta,\phi).
\end{align}
Proceeding as in the previous subsection, one writes
\begin{align}
\label{rr-4}
\BBraket{N-1}{n_1,n_2}{q^{A_0/2}B_{-}U_{BC}(\phi)U_{AB}(\theta)}{x_1,x_2}{N}=\BBraket{N-1}{n_1,n_2}{U_{BC}(q^{-1/2}\phi)U_{AB}(\theta)\;\widetilde{B}_{-}}{x_1,x_2}{N},
\end{align}
where
\begin{align}
\label{BTilde-1}
\widetilde{B}_{-}=U_{AB}^{\dagger}(\theta)\,U_{BC}^{\dagger}(q^{-1/2}\phi)\;q^{-A_0/2}B_{-}\;U_{BC}(\phi)U_{AB}(\theta).
\end{align}
Upon using \eqref{Rel1} with the $A$s replaced by $B$s and the $B$s replaced by $C$s, one finds
\begin{align*}
\widetilde{B}_{-}=U_{AB}^{\dagger}(\theta)\,\Bigg\{q^{(A_0+B_0)/2}\Big[\sqrt{1-\phi^2 q^{-1}q^{-C_0}}q^{-B_0/2}B_{-}+\phi q^{-1/2} q^{-C_0/2}C_{-}\Big]\Bigg\}\,U_{AB}(\theta).
\end{align*}
It is obvious that $q^{A_0+B_0}$ commutes with $U_{AB}(\theta)$, and similarly for all the terms involving the generators $C_0$, $C_{\pm}$. Using the conjugation formula 
\begin{align*}
U_{AB}^{\dagger}(\theta)\,q^{-B_0/2} B_{-}\,U_{AB}(\theta)=q^{-B_0/2}\,B_{-}\,\sqrt{1-\theta^2q^{-B_0}}-\theta\,q^{-B_0}A_{-} q^{-A_0/2},
\end{align*}
which can be obtained by a straightforward calculation with the help of the $q$-BCH formulas \eqref{BCH-1}, \eqref{BCH-2}, one obtains
\begin{align}
\label{rr-5}
\widetilde{B}_{-}=q^{(A_0+B_0)/2}\Big[\sqrt{1-\phi^2 q^{-1}q^{-C_0}}\left(q^{-B_0/2}\,B_{-}\,\sqrt{1-\theta^2q^{-B_0}}-\theta\,q^{-B_0}A_{-} q^{-A_0/2}\right)+\phi q^{-1/2}q^{-C_0/2}C_{-}\Big].
\end{align}
Upon inserting \eqref{rr-5} in \eqref{rr-4}, and comparing with \eqref{rr-3}, one finds
\begin{multline}
q^{n_1/2}\sqrt{1-q^{n_2+1}}\,\Xi_{n_1,n_2+1}^{(N)}(x_1,x_2;\theta,\phi)=
\phi q^{x_1+x_2}q^{-N/2}\sqrt{1-q^{N-x_1-x_2}}\,\Xi_{n_1,n_2}^{(N-1)}(x_1,x_2;\theta,q^{-1/2}\phi)
\\
-q^{-1/2}\theta q^{-x_2/2} \sqrt{1-\phi^2 q^{x_1+x_2-N-1}}\,\sqrt{1-q^{x_1}}\,\Xi_{n_1,n_2}^{(N-1)}(x_1-1,x_2;\theta,q^{-1/2}\phi)
\\
+q^{x_1/2}\sqrt{1-q^{x_2}}\sqrt{1-\theta^2 q^{-x_2}}\sqrt{1-\phi^2 q^{x_1+x_2-N-1}}\,\Xi_{n_1,n_2}^{(N-1)}(x_1,x_2-1;\theta,q^{-1/2}\phi).
\end{multline}
Using the expression \eqref{2var-Elements-Explicit}, one obtains
\begin{multline}
\label{rr-23}
q^{N-n_2}\mathbf{K}_{n_1,n_2+1}(x_1,x_2;\theta,q^{1/2}\phi;N)=(1-q^{x_1})(1-\phi^{-2}q^{N-x_1-x_2})\;\mathbf{K}_{n_1,n_2}(x_1-1,x_2;\theta,\phi;N-1)
\\
\qquad \qquad +q^{x_1}(1-q^{x_2})(1-\phi^{-2}q^{N-x_1-x_2})\;\mathbf{K}_{n_1,n_2}(x_1,x_2-1;\theta,\phi;N-1)
\\
+q^{x_1+x_2}(1-q^{N-x_1-x_2})\;\mathbf{K}_{n_1,n_2}(x_1,x_2;\theta,\phi;N-1).
\end{multline}
The relations \eqref{rr-22} and \eqref{rr-23} can be combined to generate $\mathbf{K}_{n_1,n_2}(x_1,x_2;\alpha_1,\alpha_2;N)$ recursively.
\subsubsection{Lowering relation in $n_2$}
In order to obtain a lowering relation in the second degree index $n_2$, one considers the matrix elements $\BBraket{N+1}{n_1,n_2}{B_{+}q^{A_0/2}\,U_{BC}(q^{-1/2}\phi)U_{AB}(\theta)}{x_1,x_2}{N}$. One has
\begin{align}
\BBraket{N+1}{n_1,n_2}{B_{+}q^{A_0/2}\,U_{BC}(q^{-1/2}\phi)U_{AB}(\theta)}{x_1,x_2}{N}=q^{n_1/2}\sqrt{\frac{1-q^{n_2}}{1-q}}\;\Xi_{n_1,n_2-1}^{(N)}(x_1,x_2;\theta,q^{-1/2}\phi;N),
\end{align}
and also
\begin{align}
\BBraket{N+1}{n_1,n_2}{B_{+}q^{A_0/2}\,U_{BC}(q^{-1/2}\phi)U_{AB}(\theta)}{x_1,x_2}{N}=\BBraket{N+1}{n_1,n_2}{U_{BC}(\phi)U_{AB}(\theta)\;\widetilde{B}_{+}}{x_1,x_2}{N},
\end{align}
where $\widetilde{B}_{+}$ is given by
\begin{align}
\label{Bp-Tilde}
\widetilde{B}_{+}=U_{AB}^{\dagger}(\theta)U_{BC}(\phi)\,B_{+}q^{A_0/2}\;U_{BC}(q^{-1/2}\phi) U_{AB}(\theta).
\end{align}
It is seen that \eqref{Bp-Tilde} is the complex conjugate of \eqref{BTilde-1}. Taking the complex conjugate of \eqref{rr-5}, one finds
\begin{align}
\widetilde{B}_{+}=\Big[\left(\sqrt{1-\theta^2q^{-B_0}}\,B_{+}\,q^{-B_0/2}-\theta\,q^{-A_0/2}A_{+} q^{-B_0}\right)\sqrt{1-\phi^2 q^{-1}q^{-C_0}}+\phi q^{-1/2}C_{+}q^{-C_0/2}\Big]q^{(A_0+B_0)/2}.
\end{align}
Proceeding in the same manner as before leads to the relation
\begin{multline}
q^{n_1/2}\sqrt{1-q^{n_2}}\;\Xi_{n_1,n_2-1}^{(N)}(x_1,x_2;\theta,q^{-1/2}\phi)=
\\
-\theta q^{-1/2}q^{-x_2/2}\sqrt{1-\phi^2 q^{-1}q^{x_1+x_2-N}}\sqrt{1-q^{x_1+1}}\;\Xi_{n_1,n_2}^{(N+1)}(x_1+1,x_2;\theta,\phi)
\\
\qquad \qquad +q^{x_1/2}\sqrt{1-\phi^2 q^{-1}q^{x_1+x_2-N}}\sqrt{1-\theta^2 q^{-1}q^{-x_2}}\sqrt{1-q^{x_2+1}}\;\Xi_{n_1,n_2}^{(N+1)}(x_1,x_2+1;\theta,\phi)
\\
+\phi q^{-1/2}q^{x_1+x_2}q^{-N/2}\sqrt{1-q^{N-x_1-x_2+1}}\;\Xi_{n_1,n_2}^{(N+1)}(x_1,x_2;\theta,\phi).
\end{multline}
In terms of the bivariate $q$-Krawtchouk polynomials, this relation amounts to
\begin{multline}
\frac{q^{n_1+n_2}(1-q^{n_2})}{\phi^2}\;\mathbf{K}_{n_1,n_2-1}(x_1,x_2;\theta,q^{-1/2}\phi;N)=-\theta^2 q^{x_1}\;\mathbf{K}_{n_1,n_2}(x_1+1,x_2;\theta,\phi;N+1)
\\
+\theta^2 q^{x_1}(1-\theta^{-2} q^{x_2+1})\;\mathbf{K}_{n_1,n_2}(x_1,x_2+1;\theta,\phi;N+1)+q^{x_1+x_2+1}\;\mathbf{K}_{n_1,n_2}(x_1,x_2;\theta,\phi;N+1).
\end{multline}
\subsection{Structure relations with respect to the index $x_1$}
Proceeding in a ``dual'' fashion, one can obtain structure relations in the first variable index $x_1$.
\subsubsection{Lowering relation in $x_1$}
In order to find a lowering relation in first variable index $x_1$, one considers the matrix element $\BBraket{N}{n_1,n_2}{U_{BC}(\phi)U_{AB}(\theta)\;q^{-A_0/2}A_{-}}{x_1,x_2}{N+1}$. One the one hand, one has
\begin{align}
\label{rr-8}
\BBraket{N}{n_1,n_2}{U_{BC}(\phi)U_{AB}(\theta)\;q^{-A_0/2}A_{-}}{x_1,x_2}{N+1}=q^{-(x_1-1)/2}\sqrt{\frac{1-q^{x_1}}{1-q}}
\; \Xi_{n_1,n_2}^{(N)}(x_1-1,x_2;\theta,\phi;N).
\end{align}
On the other hand, one can write
\begin{align}
\label{rr-9}
\BBraket{N}{n_1,n_2}{U_{BC}(\phi)U_{AB}(\theta)\;q^{-A_0/2}A_{-}}{x_1,x_2}{N+1}=
\BBraket{N}{n_1,n_2}{\;\widehat{A}_{-} \;U_{BC}(\phi)U_{AB}(\theta)}{x_1,x_2}{N+1},
\end{align}
where 
\begin{align}
\label{rr-10}
\widehat{A}_{-}=U_{BC}(\phi) U_{AB}(\theta)\;q^{-A_0/2}A_{-}\;U^{\dagger}_{AB}(\theta)U_{BC}^{\dagger}(\phi).
\end{align}
To compute \eqref{rr-10}, one first uses the conjugation formula
\begin{align}
\label{Conju-3}
 U_{AB}(\theta)\,q^{-A_0/2} A_{-}\,U_{AB}^{\dagger}(\theta)=q^{-A_0/2}\,A_{-}\,\sqrt{1-\theta^2q^{-A_0}}-\theta\,q^{-A_0}B_{-} q^{-B_0/2},
\end{align}
which is obtained straightforwardly from the $q$-BCH formulas \eqref{BCH-1} and \eqref{BCH-2}; this lead to
\begin{align}
\widehat{A}_{-}=U_{BC}(\phi)\left[q^{-A_0/2}\,A_{-}\,\sqrt{1-\theta^2q^{-A_0}}-\theta q^{-1/2}\,q^{-A_0}q^{-B_0/2}B_{-} \right]U_{BC}^{\dagger}(\phi).
\end{align}
Using \eqref{Conju-3} again with $A$ replaced by $B$ and $B$ replaced by $C$, one finds
\begin{align}
\label{rr-11}
\widehat{A}_{-}=\Bigg\{q^{-A_0/2}\,A_{-}\,\sqrt{1-\theta^2q^{-A_0}}-\theta q^{-1/2}\,q^{-A_0}\Big[q^{-B_0/2}\,B_{-}\,\sqrt{1-\phi^2q^{-B_0}}-\phi q^{-1/2}\,q^{-B_0}q^{-C_0/2}C_{-}\Big]\Bigg\}.
\end{align}
Upon substituting \eqref{rr-11} in \eqref{rr-9}, using the actions \eqref{Action} and comparing with \eqref{rr-8}, one finds
\begin{multline}
q^{-(x_1-1)/2}\sqrt{1-q^{x_1}}\;\Xi_{n_1,n_2}^{(N)}(x_1-1,x_2;\theta,\phi)=
\\
q^{-n_1/2}\sqrt{1-q^{n_1+1}}\sqrt{1-\theta^2 q^{-1}q^{-n_1}}\Xi_{n_1+1,n_2}^{(N+1)}(x_1,x_2;\theta,\phi)
\\
-\theta q^{-1/2}q^{-n_1}q^{-n_2/2}\sqrt{1-q^{n_2+1}}\sqrt{1-\phi^2 q^{-1} q^{-n_2}}\;\Xi_{n_1,n_2+1}^{(N+1)}(x_1,x_2;\theta,\phi)
\\
+\theta \phi q^{-1}q^{-(n_1+n_2)/2}q^{-N/2}\sqrt{1-q^{N-n_1-n_2+1}}\;\Xi_{n_1,n_2}^{(N+1)}(x_1,x_2;\theta,\phi).
\end{multline}
In terms of the polynomials, this amounts to
\begin{multline}
\frac{q(1-q^{-x_1})}{\theta^2}\;\mathbf{K}_{n_1,n_2}(x_1-1,x_2;\theta,\phi;N)=q^{-2n_1}\;\mathbf{K}_{n_1+1,n_2}(x_1,x_2;\theta,\phi;N+1)
\\+\phi^2 q^{-n_1-2n_2-1}\;\mathbf{K}_{n_1,n_2+1}(x_1,x_2;\theta,\phi;N+1)-\phi^2 q^{-N-2}(1-q^{N-n_1-n_2+1})\;\mathbf{K}_{n_1,n_2}(x_1,x_2;\theta,\phi;N+1).
\end{multline}
\subsubsection{Raising relation in $x_1$}
A raising relation in the first variable index $x_1$ can be obtained by considering the matrix element $\BBraket{N}{n_1,n_2}{U_{BC}(\phi)U_{AB}(\theta)\;A_{+}q^{-A_0/2}}{x_1,x_2}{N-1}$. One has
\begin{align}
\label{rr-12}
\BBraket{N}{n_1,n_2}{U_{BC}(\phi)U_{AB}(\theta)\;A_{+}q^{-A_0/2}}{x_1,x_2}{N-1}=q^{-x_1/2}\sqrt{\frac{1-q^{x_1+1}}{1-q}}
\; \Xi_{n_1,n_2}^{(N)}(x_1+1,x_2;\theta,\phi),
\end{align}
and also
\begin{align}
\label{rr-13}
\BBraket{N}{n_1,n_2}{U_{BC}(\phi)U_{AB}(\theta)\;A_{+}q^{-A_0/2}}{x_1,x_2}{N-1}=\BBraket{N}{n_1,n_2}{\widehat{A}_{+}\;U_{BC}(\phi)U_{AB}(\theta)}{x_1,x_2}{N+1},
\end{align}
where
\begin{align}
\label{rr-14}
\widehat{A}_{+}=U_{BC}(\phi)U_{AB}(\theta)\;A_{+}q^{-A_0/2} U_{AB}^{\dagger}(\theta) U_{BC}^{\dagger}(\phi).
\end{align}
The expression for \eqref{rr-14} is obtained directly by taking the complex conjugate of \eqref{rr-11}. Substituting the result in \eqref{rr-13} and comparing with \eqref{rr-12}, one finds
\begin{multline}
q^{-x_1/2}\sqrt{1-q^{x_1+1}} \Xi_{n_1,n_2}^{(N)}(x_1+1,x_2;\theta,\phi)=q^{1/2}q^{-n_1/2}\sqrt{1-\theta^2 q^{-n_1}}\sqrt{1-q^{n_1}}\;\Xi_{n_1-1,n_2}^{(N-1)}(x_1,x_2;\theta,\phi)
\\
-\theta q^{-n_1} q^{-n_2/2}\sqrt{1-q^{n_2}}\sqrt{1-\phi^2 q^{-n_2}}\;\Xi_{n_1,n_2-1}^{(N-1)}(x_1,x_2;\theta,\phi)
\\
+\theta \phi q^{-1/2}q^{-(N+n_1+n_2)/2}\sqrt{1-q^{N-n_1-n_2}}\;\Xi_{n_1,n_2}^{(N-1)}(x_1,x_2;\theta,\phi).
\end{multline}
Using the expression \eqref{2var-Elements-Explicit} for the matrix elements, this relations translates to
\begin{multline}
\mathbf{K}_{n_1,n_2}(x_1+1,x_2;\theta,\phi;N)=q^{-1}(1-q^{n_1})(1-\theta^{-2}q^{n_1})\;\mathbf{K}_{n_1-1,n_2}(x_1,x_2;\theta,\phi;N-1)
\\
+q^{-1-n_1}(1-q^{n_2})(1-\phi^{-2}q^{n_2})\;\mathbf{K}_{n_1,n_2-1}(x_1,x_2;\theta,\phi;N-1)+q^{-n_1-n_2}\mathbf{K}_{n_1,n_2}(x_1,x_2;\theta,\phi;N-1).
\end{multline}
\subsection{Structure relations with respect to the index $x_2$}
We now obtain structure relations with respect to the second variable index $x_2$.
\subsubsection{Lowering relation in $x_2$}
In order to obtain a lowering relation with respect to the second variable index $x_2$, consider the matrix element $\BBraket{N}{n_1,n_2}{U_{BC}(\phi)U_{AB}(\theta)\;q^{C_0/2} B_{-}}{x_1,x_2}{N+1}$. One has
\begin{align}
\label{rr-15}
\BBraket{N}{n_1,n_2}{U_{BC}(\phi)U_{AB}(\theta)\;q^{C_0/2} B_{-}}{x_1,x_2}{N+1}=q^{(N-x_1-x_2+1)/2}\sqrt{\frac{1-q^{x_2}}{1-q}}\;\Xi_{n_1,n_2}^{(N)}(x_1,x_2-1;\theta,\phi),
\end{align}
and also
\begin{align}
\label{rr-16}
\BBraket{N}{n_1,n_2}{U_{BC}(\phi)U_{AB}(\theta)\;q^{C_0/2} B_{-}}{x_1,x_2}{N+1}=\BBraket{N}{n_1,n_2}{\widehat{B}_{-}\;U_{BC}(\phi)U_{AB}(q^{1/2}\theta)}{x_1,x_2}{N+1},
\end{align}
where
\begin{align}
\label{rr-17}
\widehat{B}_{-}=U_{BC}(\phi)U_{AB}(\theta)\;q^{C_0/2}B_{-}\;U_{AB}^{\dagger}(q^{1/2}\theta) U_{BC}^{\dagger}(\phi).
\end{align}
To determine the expression for $\widehat{B}_{-}$, one needs the conjugation relation
\begin{align}
U_{AB}(\theta)\; B_{-}\; U_{AB}^{\dagger}(q^{1/2}\theta)=\sqrt{1-\theta^2q^{-A_0}}\,B_{-}+\theta q^{(B_0-A_0)/2}A_{-},
\end{align}
which is again obtained from the $q$-BCH relations. Using this relation, \eqref{rr-17} becomes
\begin{align}
\widehat{B}_{-}=U_{BC}(\phi)\;\Bigg\{q^{(B_0+C_0)/2}\Big[\sqrt{1-\theta^2q^{-A_0}}\,q^{-B_0/2}B_{-}+\theta q^{-A_0/2}A_{-}\Big]\Bigg\}\;U_{BC}^{\dagger}(\phi).
\end{align}
Upon using \eqref{Conju-3} in the above (with $A\rightarrow B$ and $B\rightarrow C$), one finds
\begin{align}
\widehat{B}_{-}=\Bigg\{q^{(B_0+C_0)/2}\Bigg\}\Bigg\{\sqrt{1-\theta^2q^{-A_0}}\,\Big[q^{-B_0/2}\,B_{-}\,\sqrt{1-\phi^2q^{-B_0}}-\phi\,q^{-B_0}C_{-} q^{-C_0/2}\Big]+\theta q^{-A_0/2}\,A_{-}\Bigg\}.
\end{align}
Substituting the above expression in \eqref{rr-16} and comparing with \eqref{rr-15}, there comes
\begin{multline}
q^{(N-x_1-x_2+1)/2}\sqrt{1-q^{x_2}}\;\Xi_{n_1,n_2}^{(N)}(x_1,x_2-1;\theta,\phi)=\theta q^{N/2}q^{-n_1}\sqrt{1-q^{n_1+1}}\;\Xi_{n_1+1,n_2}^{(N+1)}(x_1,x_2;q^{1/2}\theta,\phi)
\\
+q^{(N-n_1-n_2)/2}\sqrt{1-\theta^2 q^{-n_1}}\sqrt{1-q^{n_2+1}}\sqrt{1-\phi^2 q^{-1} q^{-n_2}}\;\Xi_{n_1,n_2+1}^{(N+1)}(x_1,x_2;q^{1/2}\theta,\phi)
\\
-\phi q^{-1/2}q^{-n_2/2}\sqrt{1-\theta^2 q^{-n_1}}\sqrt{1-q^{N-n_1-n_2+1}}\;\Xi_{n_1,n_2}^{(N+1)}(x_1,x_2;q^{1/2}\theta,\phi).
\end{multline}
In terms of the polynomials $\mathbf{K}_{n_1,n_2}(x_1,x_2;\theta,\phi;N)$ this relation is of the form
\begin{multline}
\frac{q^{1-x_1-x_2}}{\theta^2}(1-q^{x_2})\;\mathbf{K}_{n_1,n_2}(x_1,x_2-1;\theta,\phi;N)=
q^{1-2n_1}\;\mathbf{K}_{n_1+1,n_2}(x_1,x_2;q^{1/2}\theta,\phi;N+1)
\\
+\phi^2 q^{-n_1-2n_2}(1-\theta^{-2}q^{n_1})\;\mathbf{K}_{n_1,n_2+1}(x_1,x_2;q^{1/2}\theta,\phi;N+1)
\\
-\phi^2 q^{-N-1}(1-q^{N-n_1-n_2+1})(1-\theta^{-2}q^{n_1})\;\mathbf{K}_{n_1,n_2}(x_1,x_2;q^{1/2}\theta,\phi;N+1).
\end{multline}
\subsubsection{Raising relation in $x_2$}
To obtain a raising relation in the second variable index $x_2$, one considers the matrix element $\BBraket{N}{n_1,n_2}{U_{BC}(\phi)U_{AB}(q^{1/2}\theta)\; B_{+}q^{C_0/2}}{x_1,x_2}{N-1}$. On the one hand, one has 
\begin{align}
\label{rr-18}
\BBraket{N}{n_1,n_2}{U_{BC}(\phi)U_{AB}(q^{1/2}\theta)\; B_{+}q^{C_0/2}}{x_1,x_2}{N-1}=q^{(N-x_1-x_2-1)/2}\sqrt{\frac{1-q^{x_2+1}}{1-q}}\;\Xi_{n_1,n_2}^{(N)}(x_1,x_2+1;\theta,\phi),
\end{align}
and on the other hand one can write
\begin{align}
\label{rr-19}
\BBraket{N}{n_1,n_2}{U_{BC}(\phi)U_{AB}(q^{1/2}\theta)\; B_{+}q^{C_0/2}}{x_1,x_2}{N-1}=\BBraket{N}{n_1,n_2}{\;\widehat{B}_{+}\;U_{BC}(\phi) U_{AB}(\theta)}{x_1,x_2}{N-1},
\end{align}
where
\begin{align}
\widehat{B}_{+}=U_{BC}(\phi)U_{AB}(q^{1/2}\theta)\;B_{+}q^{C_0/2}\; U_{AB}^{\dagger}(\theta)U_{BC}^{\dagger}(\phi).
\end{align}
Since $\widehat{B}_{+}=\widehat{B}_{-}^{\dagger}$, one directly finds
\begin{multline}
q^{(N-x_1-x_2-1)/2}\sqrt{1-q^{x_2+1}}\;\Xi_{n_1,n_2}^{(N)}(x_1,x_2+1; q^{1/2}\theta,\phi)=
\theta q^{1/2}q^{-n_1}q^{N/2}\sqrt{1-q^{n_1}}\;\Xi_{n_1-1,n_2}^{(N-1)}(x_1,x_2;\theta,\phi)
\\
+q^{(N-n_1-n_2)/2}\sqrt{1-\phi^2 q^{-n_2}}\sqrt{1-q^{n_2}}\sqrt{1-\theta^2 q^{-n_1}}\;\Xi_{n_1,n_2-1}^{(N-1)}(x_1,x_2;\theta,\phi)
\\
-\phi q^{-1/2}q^{-n_2/2}\sqrt{1-q^{N-n_1-n_2}}\sqrt{1-\theta^2 q^{-n_1}}\;\Xi_{n_1,n_2}^{(N-1)}(x_1,x_2;\theta,\phi).
\end{multline}
In terms of the polynomials $\mathbf{K}_{n_1,n_2}(x_1,x_2;\theta,\phi;N)$, this gives
\begin{multline}
\mathbf{K}_{n_1,n_2}(x_1,x_2+1;q^{1/2}\theta,\phi;N)=q^{-1}(1-q^{n_1})\;\mathbf{K}_{n_1-1,n_2}(x_1,x_2;\theta,\phi;N-1)
\\
+q^{-n_1-1}(1-\phi^{-2}q^{n_2})(1-q^{n_2})\;\mathbf{K}_{n_1,n_2-1}(x_1,x_2;\theta,\phi;N-1)+q^{-(n_1+n_2)}\;\mathbf{K}_{n_1,n_2}(x_1,x_2;\theta,\phi;N-1).
\end{multline}
\section{Bispectrality}
In this section, the algebraic model for the two-variable $q$-Krawtchouk polynomials is used to derive their recurrence relations and their eigenvalue equations. These relations coincide with those of \cite{2011_Iliev_TransAmerMathSoc_363_1577}.
\subsection{Difference equations}
To derive the two eigenvalue equations that determine the polynomials $\mathbf{K}_{n_1,n_2}(x_1,x_2;\alpha_1,\alpha_2;N)$, one could find suitable combinations of the raising and lowering relations on the degree indices $n_1$ and $n_2$ obtained in the previous section. However, it will prove more instructive to proceed directly by finding conjugation formulas for the appropriate operators.
\subsubsection{First difference equation}
Consider the matrix element $\BBraket{N}{n_1,n_2}{q^{A_0}\;U_{BC}(\phi)U_{AB}(\theta)}{x_1,x_2}{N}$. On the one hand, one has
\begin{align}
\label{rr-24}
\BBraket{N}{n_1,n_2}{q^{A_0}\;U_{BC}(\phi)U_{AB}(\theta)}{x_1,x_2}{N}=q^{n_1}\,\Xi^{(N)}_{n_1,n_2}(x_1,x_2;\theta,\phi).
\end{align}
On the other hand, one can write
\begin{align}
\label{rr-25}
\BBraket{N}{n_1,n_2}{q^{A_0}\;U_{BC}(\phi)U_{AB}(\theta)}{x_1,x_2}{N}=
\BBraket{N}{n_1,n_2}{U_{BC}(\phi)U_{AB}(\theta)\;\widehat{q^{A_0}}}{x_1,x_2}{N},
\end{align}
where
\begin{align}
\label{rr-26}
\widehat{q^{A_0}}=U_{AB}^{\dagger}(\theta)U_{BC}^{\dagger}(\phi)\;q^{A_0}\;U_{BC}(\phi)U_{AB}(\theta)
=U_{AB}^{\dagger}(\theta)\;q^{A_0}\;U_{AB}(\theta).
\end{align}
The conjugation $U_{AB}^{\dagger}(\theta)\;q^{A_0}\;U_{AB}(\theta)$ can be performed with the help of the $q$-BCH formulas. As this computation is not straightforward, we present the details of the calculation. Using the explicit expression \eqref{UXY} for the $q$-rotation operator, one writes
\begin{multline*}
\widehat{q^{A_0}}=e_{q}^{1/2}\left(\theta^2 q^{-B_0}\right)\,e_{q}\left(\theta(1-q)q^{-(A_0+B_0)/2}A_{-}B_{+}\right)\,E_{q}\left(-\theta (1-q)q^{-(A_0+B_0)/2}A_{+}B_{-}\right)
\\
\times E_{q}^{1/2}\left(-\theta^2 q^{-A_0}\right)\;\Big\{q^{A_0}\Big\}\;e_{q}^{1/2}\left(\theta^2 q^{-A_0}\right)
\\
\times e_{q}\left(\theta(1-q)q^{-(A_0+B_0)/2}A_{+}B_{-}\right)\,E_{q}\left(-\theta(1-q)q^{-(A_0+B_0)/2}A_{-}B_{+}\right)\,E_{q}^{1/2}\left(-\theta^2 q^{-B_0}\right).
\end{multline*}
The first conjugation is trivial and one can immediately write
\begin{multline*}
\widehat{q^{A_0}}=e_{q}^{1/2}\left(\theta^2 q^{-B_0}\right)\,e_{q}\left(\theta(1-q)q^{-(A_0+B_0)/2}A_{-}B_{+}\right)\,
\\
\times E_{q}\left(-\theta (1-q)q^{-(A_0+B_0)/2}A_{+}B_{-}\right)\;q^{A_0}\;e_{q}\left(\theta(1-q)q^{-(A_0+B_0)/2}A_{+}B_{-}\right)
\\
\times E_{q}\left(-\theta(1-q)q^{-(A_0+B_0)/2}A_{-}B_{+}\right)\,E_{q}^{1/2}\left(-\theta^2 q^{-B_0}\right).
\end{multline*}
The next step is performed using the identity $q^{A_0}\,e_{q}(\lambda\,A_{+})=e_{q}(q\lambda A_{+})\,q^{A_0}$; one then finds
\begin{align*}
\widehat{q^{A_0}}=e_{q}^{1/2}\left(\theta^2 q^{-B_0}\right)\;J_1\,J_2\;E_{q}^{1/2}\left(-\theta^2 q^{-B_0}\right),
\end{align*}
where
\begin{align}
\label{J1}
J_1=e_{q}\left(\theta(1-q)q^{-(A_0+B_0)/2}A_{-}B_{+}\right) \Bigg[1-\theta(1-q)q^{-(A_0+B_0)/2}A_{+}B_{-}\Bigg]E_{q}\left(-\theta(1-q)q^{-(A_0+B_0)/2}A_{-}B_{+}\right),
\end{align}
and 
\begin{align}
\label{J2}
J_2=e_{q}\left(\theta(1-q)q^{-(A_0+B_0)/2}A_{-}B_{+}\right) \Bigg[q^{A_0}\Bigg]E_{q}\left(-\theta(1-q)q^{-(A_0+B_0)/2}A_{-}B_{+}\right).
\end{align}
The expression for \eqref{J1} can be obtained from the conjugation formula
\begin{align*}
e_{q}(\alpha\,A_{-}B_{+})\,A_{+}B_{-}\,E_{q}(-\alpha\,A_{-}B_{+})=A_{+}B_{-}+\frac{\alpha}{(1-q)^2}\,q^{A_0}\,-\frac{\alpha}{(1-q)^2}\frac{1}{1-\alpha\,A_-B_{+}}q^{B_0},
\end{align*}
which follows straightforwardly from \eqref{BCH-2}. One thus has
\begin{multline*}
J_1=1-\theta\,(1-q)\,q^{-(A_0+B_0)/2}\,\Bigg[A_{+}B_{-}+\frac{\theta\,q^{-(A_0+B_0)/2}}{(1-q)}\,q^{A_0}
\\
-\frac{\theta\,q^{-(A_0+B_0)/2}}{(1-q)}\,\frac{1}{1-\theta\,(1-q)\,q^{-(A_0+B_0)/2}A_{-}B_{+}}\,q^{B_0}\Bigg].
\end{multline*}
One easily finds
\begin{align*}
J_2=q^{A_0}\left[1-\theta\,(1-q)\,q^{-(A_0+B_0)/2}A_{-}B_{+}\right].
\end{align*}
Performing simplifications on the product of $J_1$ and $J_2$ using the commutation relations \eqref{qOsc-Algebra}, one finds
\begin{multline*}
\widehat{q^{A_0}}=e_{q}^{1/2}\left(\theta^2 q^{-B_0}\right)\;\Bigg\{q^{A_0}[1-\theta^2q^{-B_0}]-\theta(1-q)q^{(A_0-B_0)/2}[1-\theta^2 q^{-B_0}]A_{-}B_{+}
\\
-\theta (1-q) A_{+}B_{-}q^{(A_0-B_0)/2}
 +q^{A_0}\left(\theta^2-q^{-1}\theta^2\,q^{-B_0}+\theta^2 q^{-1}q^{-(A_0+B_0)}\right)\Bigg\}\;E_{q}^{1/2}\left(-\theta^2 q^{-B_0}\right).
\end{multline*}
The last conjugation is easily carried out; the final result is 
\begin{multline}
\label{Conju-A0}
\widehat{q^{A_0}}=q^{A_0}\left[1-(1+q^{-1})\theta^2 q^{-B_0}+\theta^2+q^{-1}\theta^2 q^{-(A_0+B_0)}\right]
 \\
 -\theta(1-q)q^{(A_0-B_0)/2}\sqrt{1-\theta^2 q^{-B_0}}A_{-}B_{+}-\theta (1-q) A_{+}B_{-}\sqrt{1-\theta^2 q^{-B_0}}q^{(A_0-B_0)/2}.
\end{multline}
Upon using \eqref{Conju-A0} in \eqref{rr-25}, using the actions \eqref{Action} and comparing with \eqref{rr-24}, one finds that the matrix elements $\Xi_{n_1,n_2}^{(N)}(x_1,x_2;\theta,\phi)=\Xi_{n_1,n_2}^{(N)}(x_1,x_2)$ satisfy the eigenvalue equation
\begin{multline}
\label{Eigen-1}
q^{n_1}\Xi_{n_1,n_2}^{(N)}(x_1,x_2)=\left[q^{x_1}(1+\theta^2)+q^{-1}\theta^2 q^{-x_2}-(1+q^{-1})\theta^2 q^{x_1-x_2}\right]\,\Xi_{n_1,n_2}^{(N)}(x_1,x_2)
 \\
 -\theta\,q^{-1} q^{(x_1-x_2)/2}\left[(1-q^{x_2+1})(1-q^{x_1})(1-\theta^2 q^{-(x_2+1)})\right]^{1/2}\,\Xi_{n_1,n_2}^{(N)}(x_1-1,x_2+1)
 \\
 -\theta q^{(x_1-x_2)/2}\left[(1-\theta^2 q^{-x_2})(1-q^{x_2})(1-q^{x_1+1})\right]^{1/2}\,\Xi_{n_1,n_2}^{(N)}(x_1+1,x_2-1).
\end{multline}
Upon using the expression \eqref{2var-Elements-Explicit} for the matrix elements, one finds that the two-variable quantum $q$-Krawtchouk polynomials $\mathbf{K}_{n_1,n_2}(x_1,x_2;\theta,\phi;N)$ satisfy the 3-point eigenvalue equation
\begin{multline}
\label{Diff-1}
 q^{n_1}\;\mathbf{K}_{n_1,n_2}(x_1,x_2)
   =\left[q^{x_1}(1+\theta^2)+q^{-1}\theta^2 q^{-x_2}-(1+q^{-1})\theta^2 q^{x_1-x_2}\right]\;\mathbf{K}_{n_1,n_2}(x_1,x_2)
 \\
 +(1-\theta^2 q^{-x_2-1})(1-q^{x_1})\;\mathbf{K}_{n_1,n_2}(x_1-1,x_2+1)
 +\theta^2\,q^{x_1-x_2}(1-q^{x_2})\;\mathbf{K}_{n_1,n_2}(x_1+1,x_2-1).
\end{multline}
\subsubsection{Second difference equation}
Another difference equation for the matrix elements $\Xi_{n_1,n_2}^{(N)}(x_1,x_2)$ and the two-variable Krawtchouk polynomials can be obtained from the matrix element $\BBraket{N}{n_1,n_2}{q^{A_0+B_0}\;U_{BC}(\phi)U_{AB}(\theta)}{x_1,x_2}{N}$. One can write this matrix elements as
\begin{align}
\BBraket{N}{n_1,n_2}{q^{A_0+B_0}\;U_{BC}(\phi)U_{AB}(\theta)}{x_1,x_2}{N}=q^{n_1+n_2}\Xi_{n_1,n_2}^{(N)}(x_1,x_2;\theta,\phi).
\end{align}
but also as
\begin{align}
\BBraket{N}{n_1,n_2}{q^{A_0+B_0}\;U_{BC}(\phi)U_{AB}(\theta)}{x_1,x_2}{N}=
\BBraket{N}{n_1,n_2}{U_{BC}(\phi)U_{AB}(\theta)\;\widehat{q^{A_0+B_0}}}{x_1,x_2}{N},
\end{align}
where
\begin{align}
\widehat{q^{A_0+B_0}}=U_{AB}^{\dagger}(\theta)U_{BC}^{\dagger}(\phi)\;q^{A_0+B_0}\;U_{BC}(\phi) U_{AB}(\theta).
\end{align}
Using the result \eqref{Conju-A0} with the $A$s replaced by $B$s and the $B$s replaced by $C$s, one has
\begin{multline}
\widehat{q^{A_0+B_0}}=U_{AB}^{\dagger}(\theta)\;q^{A_0+B_0}\,\Bigg\{\left[1-(1+q^{-1})\phi^2 q^{-C_0}+\phi^2+q^{-1}\phi^2 q^{-(B_0+C_0)}\right]
 \\
 -\phi(1-q)q^{-(B_0+C_0)/2}\sqrt{1-\phi^2 q^{-C_0}}B_{-}C_{+}-\phi\,q^{-1}\,(1-q) B_{+}C_{-}\,\sqrt{1-\phi^2 q^{-C_0}}q^{-(B_0+C_0)/2}\Bigg\}\,U_{AB}(\theta).
\end{multline}
To complete this calculation, one needs to use the conjugation formula
\begin{multline*}
U_{AB}^{\dagger}(\theta)\,q^{-B_0}\,U_{AB}(\theta)=q^{-B_0}\Bigg\{\left[1-(1+q^{-1})\theta^2 q^{-B_0}+\theta^2+q^{-1}\theta^2 q^{-(A_0+B_0)}\right]
 \\
 -\theta(1-q)q^{-(A_0+B_0)/2}\sqrt{1-\theta^2 q^{-B_0}}A_{-}B_{+}-\theta\,q^{-1}\,(1-q) A_{+}B_{-}\sqrt{1-\theta^2 q^{-B_0}}q^{-(A_0+B_0)/2}\Bigg\},
\end{multline*}
which can be obtained in a fashion similar to the way \eqref{rr-26} was derived, as well as the formula
\begin{align*}
U_{AB}^{\dagger}(\theta)\,q^{-B_0/2} B_{-}\,U_{AB}(\theta)=q^{-B_0/2}\,B_{-}\,\sqrt{1-\theta^2q^{-B_0}}-\theta\,q^{-B_0}A_{-} q^{-A_0/2},
\end{align*}
and its complex conjugate, which are easily obtained from the $q$-BCH relations. The final result is
\begin{multline}
\widehat{q^{A_0+B_0}}=q^{A_0+B_0}\left\{1-(1+q^{-1})\phi^2 q^{-C_0}+\phi^2\right\}
\\
+q^{-1}\phi^2 q^{A_0-C_0}\left[1-(1+q^{-1})\theta^2 q^{-B_0}+\theta^2+q^{-1}\theta^2 q^{-(A_0+B_0)}\right]
  \\
  -q^{-1}\phi^2\theta(1-q)q^{(A_0-B_0-2C_0)/2}\sqrt{1-\theta^2 q^{-B_0}}A_{-}B_{+}-q^{-1}\phi^2\theta(1-q)A_{+}B_{-}\sqrt{1-\theta^2 q^{-B_0}}q^{(A_0-B_0-2C_0)/2}
  \\
  -\phi(1-q)q^{(2A_0+B_0)/2}B_{-}\sqrt{1-\theta^2 q^{-B_0}}q^{-C_0/2}\sqrt{1-\phi^2 q^{-C_0}}C_{+}
  \\
  -\phi(1-q)\sqrt{1-\theta^2 q^{-B_0}}B_{+}q^{(2A_0+B_0)/2}C_{-}\sqrt{1-\phi^2 q^{-C_0}}q^{-C_0/2}
  \\
  +\theta\,\phi\,q^{-1}(1-q)A_{-}q^{A_0/2} q^{-C_0/2}\sqrt{1-\phi^2 q^{-C_0}}C_{+}
  +\theta\,\phi\,q^{-1}(1-q)q^{A_0/2}A_{+}C_{-}\sqrt{1-\phi^2 q^{-C_0}} q^{-C_0/2}.
\end{multline}
As a consequence, we have the following 7-point eigenvalue equation
\begin{multline}
\label{Eigen-2}
 q^{n_1+n_2}\,\Xi_{n_1,n_2}^{(N)}(x_1,x_2)=
 \bigg[q^{x_1+x_2}\left(1-(1+q^{-1})\phi^2 q^{x_1+x_2-N}+\phi^2\right)
 \\
 +q^{-1}\phi^2q^{2x_1+x_2-N}\left(1-(1+q^{-1})\theta^2q^{-x_2}+\theta^2 +q^{-1}\theta^2 q^{-(x_1+x_2)}\right)\bigg]\;\Xi_{n_1,n_2}^{(N)}(x_1,x_2)
 \\
 -q^{-2}\theta\,\phi^2q^{(3x_1+x_2-2N)/2}\sqrt{(1-q^{x_2+1})(1-q^{x_1})(1-\theta^2 q^{-(x_2+1)})}\;\Xi_{n_1,n_2}^{(N)}(x_1-1,x_2+1)
 \\
 -q^{-1}\,\theta\,\phi^2 q^{(3x_1+x_2-2N)/2}\sqrt{(1-q^{x_2})(1-q^{x_1+1})(1-\theta^2 q^{-x_2})}\;\Xi_{n_1,n_2}^{(N)}(x_1+1,x_2-1)
 \\
 -\phi\,q^{(3x_1+2x_2-N-2)/2}\sqrt{(1-q^{x_2})(1-\theta^2 q^{-x_2})(1-q^{N-x_1-x_2+1})(1-\phi^2 q^{x_1+x_2-N-1})}\;\Xi_{n_1,n_2}^{(N)}(x_1,x_2-1)
 \\
 -\phi\,q^{(3x_1+2x_2-N)/2}\sqrt{(1-q^{x_2+1})(1-\theta^2 q^{-(x_2+1)})(1-q^{N-x_1-x_2})(1-\phi^2 q^{x_1+x_2-N})}\;\Xi_{n_1,n_2}^{(N)}(x_1,x_2+1)
 \\
+q^{-1}\theta\,\phi q^{(2x_1+x_2-N-1)/2}\sqrt{(1-q^{x_1})(1-q^{N-x_1-x_2+1})(1-\phi^2 q^{x_1+x_2-N-1})}\;\Xi_{n_1,n_2}^{(N)}(x_1-1,x_2)
 \\
 +q^{-1}\theta\,\phi q^{(2x_1+x_2-N+1)/2}\sqrt{(1-q^{x_1+1})(1-q^{N-x_1-x_2})(1-\phi^2 q^{x_1+x_2-N})}\;\Xi_{n_1,n_2}^{(N)}(x_1+1,x_2).
\end{multline}
In terms of the polynomials, this relation amounts to
\begin{multline}
\label{Diff-2}
q^{n_1+n_2}\,\mathbf{K}_{n_1,n_2}(x_1,x_2)=
 \bigg[q^{x_1+x_2}\left(1-(1+q^{-1})\phi^2 q^{x_1+x_2-N}+\phi^2\right)
 \\
 +q^{-1}\phi^2q^{2x_1+x_2-N}\left(1-(1+q^{-1})\theta^2q^{-x_2}+\theta^2 +q^{-1}\theta^2 q^{-(x_1+x_2)}\right)\bigg]\;\mathbf{K}_{n_1,n_2}(x_1,x_2)
 \\
 +\phi^2 q^{x_1+x_2-N-1}(1-q^{x_1})(1-\theta^2 q^{-x_2-1})\;\mathbf{K}_{n_1,n_2}(x_1-1,x_2+1)
 \\
 \hfill +\theta^2 \phi^2 q^{2x_1-N-1}(1-q^{x_2})\;\mathbf{K}_{n_1,n_2}(x_1+1,x_2-1)
 \\
 +q^{x_1}(1-q^{x_2})(1-\phi^2 q^{x_1+x_2-N-1})\; \mathbf{K}_{n_1,n_2}(x_1,x_2-1)
 \\
 \hfill-\phi^2 q^{x_1+x_2}(1-q^{x_1+x_2-N})(1-\theta^2 q^{-x_2-1})\; \mathbf{K}_{n_1,n_2}(x_1,x_2+1)
 \\
+(1-q^{x_1})(1-\phi^2 q^{x_1+x_2-N-1})\;\mathbf{K}_{n_1,n_2}(x_1-1,x_2)- \theta^2\phi^2 q^{x_1-1}(1-q^{x_1+x_2-N})\;\mathbf{K}_{n_1,n_2}(x_1+1,x_2).
\end{multline}
The eigenvalue equations \eqref{Diff-1} and \eqref{Diff-2} coincide with those of \cite{2011_Iliev_TransAmerMathSoc_363_1577}.
\subsection{Recurrence relations}
We shall now derive recurrence relations for the matrix elements $\Xi_{n_1,n_2}^{(N)}(x_1,x_2)$ as well as for the two-variable polynomials $\mathbf{K}_{n_1,n_2}(x_1,x_2)$.
\subsubsection{First recurrence relation}
To obtain a first recurrence relation for the matrix elements $\Xi_{n_1,n_2}^{(N)}(x_1,x_2)$, one can start from the eigenvalue equation \eqref{Eigen-1} and apply the duality relation \eqref{duality-2} to exchange the roles played by the degree indices $n_1$, $n_2$ and the variable indices $x_1$, $x_2$. Applying \eqref{duality-2} on \eqref{Eigen-1} and performing the change of variables and parameters
\begin{align*}
x_2\leftrightarrow n_2,\qquad N-n_1-n_2\leftrightarrow N-x_1-x_2,\qquad \theta\leftrightarrow \phi,
\end{align*}
one finds
\begin{multline*}
q^{-x_1-x_2}\,\Xi_{n_1,n_2}^{(N)}(x_1,x_2)=-\phi q^{-(n_1+2n_2-N)/2} \sqrt{(1-q^{N-n_1-n_2+1})(1-q^{n_2})(1-\phi^2 q^{-n_2})}\;\Xi_{n_1,n_2-1}^{(N)}(x_1,x_2)
\\
+q^{-(n_1+2n_2+1)}\left[q^{n_2+1}(1+\phi^2)+\phi^2 q^{n_1+n_2-N}-(1+q)\;\phi^2\right]\;\Xi_{n_1,n_2}^{(N)}(x_1,x_2)
\\
-\phi q^{-(n_1+2n_2+N-2)/2}\sqrt{(1-q^{n_2+1})(1-q^{N-n_1-n_2})(1-\phi^2 q^{-n_2-1})}\;\Xi_{n_1,n_2+1}^{(N)}(x_1,x_2).
\end{multline*}
In terms of the polynomials, one finds
\begin{multline}
q^{-x_1-x_2}\;\mathbf{K}_{n_1,n_2}(x_1,x_2)=
   q^{-(n_1+2n_1+1)}\Big[q^{n_2+1}(1+\phi^2)+\phi^2 q^{n_1+n_2-N}-(1+q)\;\phi^2\Big]\;\mathbf{K}_{n_1,n_2}(x_1,x_2)
   \\
   +q^{-(n_1+1)}(1-q^{n_2})(1-q^{n_1+n_2-N-1})(1-\phi^2 q^{-n_2})\,\mathbf{K}_{n_1,n_2-1}(x_1,x_2)
   -q^{-(3n_2+n_1+1)}\phi^2 \,\mathbf{K}_{n_1,n_2+1}(x_1,x_2).
\end{multline}
\subsubsection{Second recurrence relation}
To obtain the second recurrence relation, one applies the same procedure as above on the eigenvalue equation \eqref{Eigen-2} with the result
\begin{multline}
q^{-x_1}\;\Xi_{n_1,n_2}^{(N)}(x_1,x_2)=q^{-(2n_1+n_2+N+1)/2}\theta \phi \sqrt{(1-q^{n_1})(1-q^{N-n_1-n_2+1})(1-\theta^2 q^{-n_1})}\;\Xi_{n_1-1,n_2}^{(N)}(x_1,x_2)
\\
-q^{-(3n_1+n_2)/2}\theta \sqrt{(1-q^{n_2+1})(1-q^{n_1})(1-\theta^2 q^{-n_1})(1-\phi^2 q^{-n_2-1})}\;\Xi_{n_1-1,n_2+1}^{(N)}(x_1,x_2)
\\
-q^{-1/2}q^{-(2n_1+n_2)}\theta^2 \phi \sqrt{(1-q^{-n_2})(1-q^{n_1+n_2-N-1})(1-\phi^2 q^{-n_2})}\;\Xi_{n_1,n_2-1}^{(N)}(x_1,x_2)
\\
+q^{-2n_1}\Bigg[q^{n_1}(1+\theta^2 -q^{-1-n_1}(1+q)\theta^2)+q^{-1-n_2}\theta^2(1+\phi^2+q^{n_1-N-1}\phi^2-q^{-n_2-1}(1+q)\phi^2)\Bigg]\;\Xi_{n_1,n_2}^{(N)}(x_1,x_2)
\\
-q^{-(3n_1+4n_2+N+4)/2}\theta^2 \phi \sqrt{(1-q^{N-n_1-n_2})(1-q^{n_2+1})(1-\phi^2 q^{-n_2-1})}\;\Xi_{n_1,n_2+1}^{(N)}(x_1,x_2)
\\
-q^{-(3n_1+n_2+1)/2} \theta \sqrt{(1-q^{n_2})(1-q^{n_1+1})(1-\theta^2 q^{-n_1-1})(1-\phi^2 q^{-n_2})}\;\Xi_{n_1+1,n_2-1}^{(N)}(x_1,x_2)
\\
+q^{(2n_1+n_2+N+3)/2} \theta \phi \sqrt{(1-q^{n_1+1})(1-q^{N-n_1-n_2})(1-\theta^2 q^{-n_1-1})}\;\Xi_{n_1+1,n_2}^{(N)}(x_1,x_2).
\end{multline}
In terms of the polynomials, one finds 
\begin{multline}
q^{-x_1}\;\mathbf{K}_{n_1,n_2}(x_1,x_2)=q^{-2n_1}\Bigg\{q^{n_1}(1+\theta^2-q^{-n_1-1}(1+q)\theta^2)
   \\
  +q^{-n_2-1}\theta^2(1+\phi^2+q^{n_1-N-1}\phi^2-q^{-n_2-1}(1+q)\phi^2) \Bigg\}\;\mathbf{K}_{n_1,n_2}(x_1,x_2)
  \\
  \hfill
  +q^{-(n_2+2)}\phi^2 (1-q^{n_1})(1-q^{n_1+n_2-N-1})(1-\theta^2 q^{-n_1})\;\mathbf{K}_{n_1-1,n_2}(x_1,x_2)
  \\
  +q^{-2(n_2+1)}\phi^2 (1-q^{n_1})(1-\theta^2 q^{-n_1})\,\mathbf{K}_{n_1-1,n_2+1}(x_1,x_2)
  \\
  \hfill
  +q^{-(2n_1+2)}\theta^2 (1-q^{n_2})(1-q^{n_1+n_2-N-1})(1-\phi^2 q^{-n_2})\,\mathbf{K}_{n_1,n_2-1}(x_1,x_2)
  \\
  -q^{-(3n_2+2n_1+2)}\theta^{2}\phi^2\; \mathbf{K}_{n_1,n_2+1}(x_1,x_2) -q^{-(3n_1+n_2+1)}\theta^2\;\mathbf{K}_{n_1+1,n_2}(x_1,x_2)
  \\
  \hfill
  +\theta^2\phi^{-2}\,q^{n_2-3n_1-2}(1-q^{n_2})(1-\phi^2 q^{-n_2})\;\mathbf{K}_{n_1+1,n_2-1}(x_1,x_2).
\end{multline}
\textbf{Remark.} Let us note that the structure relations with respect to the variable indices $x_1,x_2$ derived in Subsections (4.3) and (4.4) cannot be obtained by the same duality argument from the structure relations with respect to the degree indices $n_1,n_2$ given in Subsections (4.1) and (4.2).
\section{Conclusion}
In this paper, we have constructed an algebraic model for the two-variable quantum $q$-Krawtchouk polynomials. We have shown that these polynomials arise as matrix elements of a three-dimensional $q$-rotation operator expressed in terms of $q$-exponentials that act on the states of three independent $q$-oscillators. We have explained how this model can be used to derive the main properties of the polynomials: orthogonality relation, duality property, structure relations, eigenvalue equations and recurrence relations.

It is apparent that the approach presented here can be extended to any number variables; the only essential difficulty lies in the complicated notation required to describe multivariate polynomials. Nevertheless, let us exhibit how the multivariate quantum $q$-Krawtchouk polynomials of Gasper and Rahman arise in the present framework. Consider $d+1$ independent copies of the $q$-oscillator algebra with generators $A_{0}^{(i)}$, $A_{\pm}^{(i)}$ and commutation relations
\begin{align*}
[A_{0}^{(i)}, A_{\pm}^{(j)}]=\pm \delta_{ij}\;A_{\pm}^{(i)},\qquad [A_{-}^{(i)},A_{+}^{(i)}]=\delta_{ij}\;q^{A_0^{(i)}},\qquad i=1,2,\ldots, d+1.
\end{align*}
Let $\{n_1,n_2\ldots,n_{d}\}$ and $\{x_1,x_2,\ldots, x_{d}\}$ be two sets of non-negative integers such that $\sum_{k=1}^{d}n_{d}\leq N$ and $\sum_{k=1}^{d}x_{k}\leq N$, where $N$ is a positive integer. We define $\mathbf{n}=(n_1,n_2,\ldots,n_{d+1})$ with $n_{d+1}=N-\sum_{k=1}^{d}n_{d}$, and similarly for $\mathbf{x}$. We shall use the notation:
\begin{align*}
|\mathbf{y}_{k}|=y_1+y_2+\cdots y_{k}.
\end{align*}
Note that by definition, we have $|\mathbf{y}_{d+1}|=N$ and $|\mathbf{y}_0|=0$. Consider the representation of the algebra generated by the $d+1$ independent $q$-oscillators on the states
\begin{align*}
\kket{\mathbf{n}}{N}=\bigotimes_{k=1}^{d+1}\ket{n_k},
\end{align*}
and defined by the following actions of the generators on the factors of the direct product:
\begin{align}
A_{+}^{(i)}\ket{n_{i}}&=\sqrt{\frac{1-q^{n_{i}+1}}{1-q}} \ket{n_{i}+1},\qquad A_{-}^{(i)}\ket{n_{i}}=\sqrt{\frac{1-q^{n_{i}}}{1-q}}\ket{n_{i}-1},\qquad 
A_0^{(i)}\ket{n_{i}}=n_{i}\ket{n_{i}}.
\end{align}
We introduce the $q$-rotation operator in the $(i,j)$ plane which is defined as
\begin{multline*}
U_{i,j}(\theta)=e_{q}^{1/2}\left(\theta^2 q^{-A_0^{(i)}}\right)\,e_{q}\left(\theta(1-q)q^{-(A_0^{(i)}+A_0^{(j)})/2}A_{+}^{(i)}A_{-}^{(j)}\right)
\\\times 
E_{q}\left(-\theta(1-q)q^{-(A_0^{(i)}+A_0^{(j)})/2}A_{-}^{(i)}A_{+}^{(j)}\right)\,E_{q}^{1/2}\left(-\theta^2 q^{-A_0^{(j)}}\right).
\end{multline*}
Let $\Xi_{\mathbf{n}}(\mathbf{x};\theta_{1},\cdots, \theta_{d})$ be the following matrix elements:
\begin{align*}
\Xi_{\mathbf{n}}(\mathbf{x};\theta_{1},\cdots, \theta_{d})
=\BBraket{N}{\mathbf{n}}{\;\overrightarrow{\prod}_{k=0}^{d-1} \;U_{d-k,d-k+1}(\theta_{d-k})\;}{\mathbf{x}}{N}.
\end{align*}
A direct calculation shows that in terms of the one-variable matrix elements  \eqref{Univariate-Elements}, the matrix elements $\Xi_{\mathbf{n}}(\mathbf{x};\theta_{1},\cdots, \theta_{d})$ have the expression
\begin{align*}
\Xi_{\mathbf{n}}(\mathbf{x};\theta_{1},\cdots, \theta_{d})=
\prod_{k=1}^{d}\xi_{n_{k}, |\mathbf{x}_{k}|-|\mathbf{n}_{k-1}|}^{(|\mathbf{x}_{k+1}|-|\mathbf{n}_{k-1}|)}(\theta_{k}).
\end{align*}
In view of \eqref{1var-Explicit}, this means that these matrix elements are proportional to the multivariate quantum $q$-Krawtchouk polynomials $\mathbf{K}_{\mathbf{n}}(\mathbf{x}; \theta_1,\ldots,\theta_{d};N)$ which read
\begin{align*}
\mathbf{K}_{\mathbf{n}}(\mathbf{x}; \theta_1,\ldots,\theta_{d};N)=
\prod_{k=1}^{d}
k_{n_{k}}\left(|\mathbf{x}_{k}|-|\mathbf{n}_{k-1}|, \frac{1}{\theta_k^2},|\mathbf{x}_{k+1}|-|\mathbf{n}_{k-1}|;q \right),
\end{align*}
where $k_{n}(x,p,N;q)$ are the one-variable quantum $q$-Krawtchouk polynomials given by \eqref{1var-qKraw}. With the above expression for $\Xi_{\mathbf{n}}(\mathbf{x};\theta_{1},\cdots, \theta_{d})$, one can proceed with the determination of the properties of the multivariate $q$-Krawtchouk polynomials.

As explained previously, the results presented here can be viewed as a $q$-generalization of \cite{2013_Genest&Vinet&Zhedanov_JPhysA_46_505203} (in the special Tratnik case), where an algebraic interpretation of Griffiths' multivariate Krawtchouk polynomials was given in terms of the rotation group acting on oscillator states. The approach of \cite{2013_Genest&Vinet&Zhedanov_JPhysA_46_505203} was generalized  in \cite{2014_Genest&Miki&Vinet&Zhedanov_JPhysA_47_215204}  and  \cite{2014_Genest&Miki&Vinet&Zhedanov_JPhysA_47_045207}, where algebraic interpretations of the multivariate Charlier and Meixner polynomials involving the Euclidean and pseudo-rotation groups were found. In light of these results, it would be of interest to generalize the approach developed here to also obtain algebraic models for the multivariate $q$-Charlier and $q$-Meixner polynomials.

\section*{Acknowledgments}
The authors would like to thank A. Zhedanov for stimulating discussions. VXG is supported by a postdoctoral fellowship from the Natural Sciences and Engineering Research Council of Canada (NSERC). The research of LV is supported in part by NSERC.
\small

\end{document}